\documentclass[leqno]{elsarticle}
\usepackage[utf8x]{inputenc}
\usepackage{amssymb}
\usepackage{amsmath}
\usepackage{newlfont}
\usepackage{graphicx}
\usepackage{color}  
\usepackage{xspace}
\usepackage{subcaption}
\usepackage{siunitx}
\usepackage{bm}
\usepackage{multirow}
\usepackage[numbers]{natbib}
\usepackage{tikz}
\usepackage[outdir=./]{epstopdf}
\usepackage{stmaryrd}
\usepackage{esint}
\usepackage[english]{babel}
\usepackage{mathtools}
\usepackage{float}
\usepackage{marginnote}
\usepackage{tabu}
\usepackage{placeins}
\usepackage{hyperref}
\usepackage{bookmark}
\usepackage{cleveref}
\usepackage{accents}
\usepackage{mathrsfs}

\usetikzlibrary{decorations.pathreplacing}
\usetikzlibrary{fadings}

\let\cite=\citet
\DeclarePairedDelimiter\abs{\lvert}{\rvert}

\usepackage{./mydef}

\newcommand{\revone}[1]{\textcolor{black}{#1}}
\newcommand{\revtwo}[1]{\textcolor{black}{#1}}
\newcommand{\modif}[1]{\textcolor{black}{#1}}

\begin{document}

	\title{Electroconvection of Thin Liquid Crystals: Model Reduction and Numerical Simulations}
	\tnotetext[t1]{A.B. is partially supported by NSF grant DMS-1817691. 
		P.W.  is partially supported by DMS-1254618, DMS-1817691 and the TAMU Hagler Institute for Advanced Study HEEP Graduate Fellowship.}
	
	\author[1]%
	{Andrea Bonito\corref{cor1}}%
	\ead{bonito@math.tamu.edu}
	
	\author[1]%
	{Peng Wei}%
	\ead{weip@math.tamu.edu}

	\cortext[cor1]{Corresponding author}
	\address[1]%
	{Department of Mathematics, Texas A\&M University, College Station, TX 77843, USA; \{bonito,weip\}@math.tamu.edu}%
	
	\begin{abstract}
		We propose a finite element method for the numerical simulation of electroconvection of thin liquid crystals.
		The liquid is located in between two concentric circular electrodes which are either assumed to be of infinite height or slim. 
		Each configuration  results in a different nonlocal electro-magnetic model defined on a two dimensional bounded domain. 
		The numerical method consists in approximating the surface charge density, the liquid velocity and pressure, and the electric potential in the two dimensional liquid region.
		Finite elements for the space discretization coupled with standard time stepping methods are put forward.
		Unlike for the infinite electrodes configuration, our numerical simulations indicate that slim electrodes are favorable for electroconvection to occur and are able to sustain the phenomena over long period of time.
		Furthermore, we provide a numerical study on the influence of the three main parameters of the system: the Rayleigh number, the Prandtl number and the electrodes aspect ratio.
	\end{abstract}
	

	\maketitle

\section{Introduction}
	The electrical convection, or \textit{electroconvection} in short, 
    is the appearance of pairs of vortices caused by an electric field  across a thin-layer of electrically charged fluid.
    	We refer to \revtwo{~\cite{morris1991rectangular,daya1997rectangular,deyirmenjian1997rectangular,deyirmenjian2005annular,daya1998annularshear,daya1999annularshear,daya2002annularshear,tsai2004annularshear,tsai2008ec}} for experimental studies in different geometries
 and to \cite{KSH2010application,KPLH2013application} for applications of electroconvection in bio-technologies.

\begin{figure} [ht]
		\begin{center}
			\begin{tabular}{cc}
				\includegraphics[scale=.12]
				{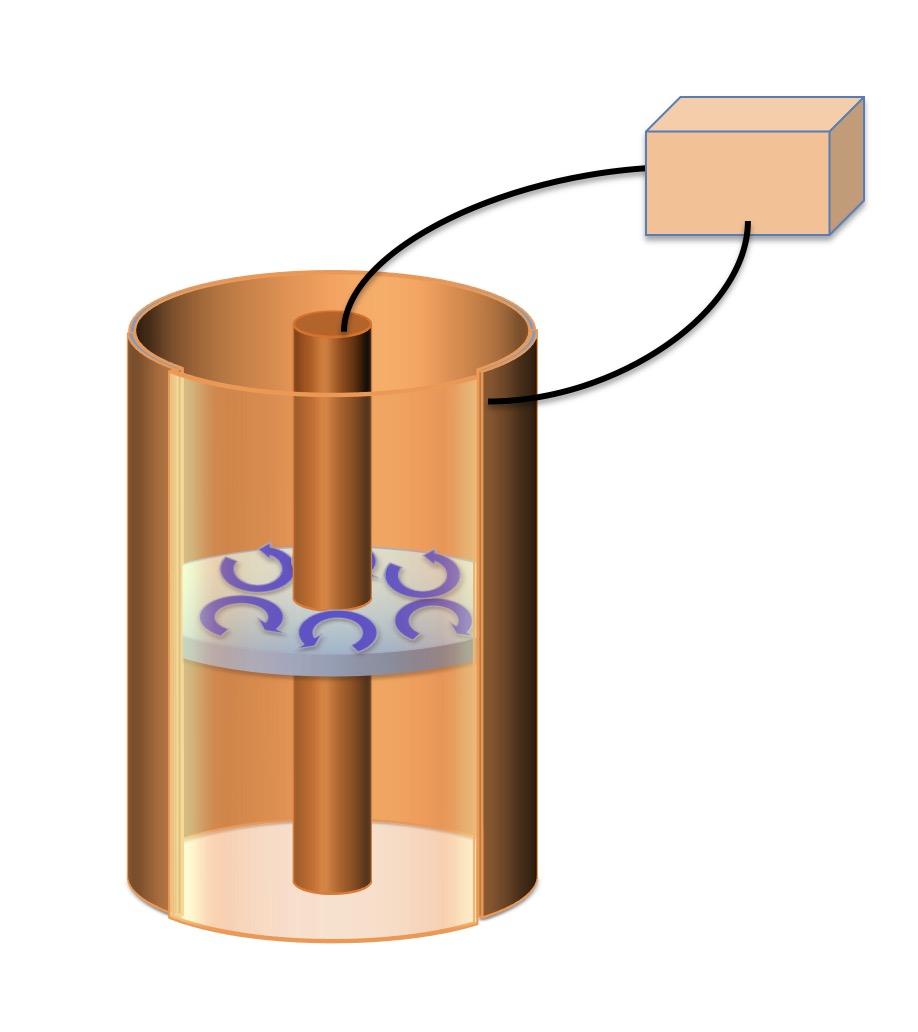}&
				\includegraphics[scale=.15]
				{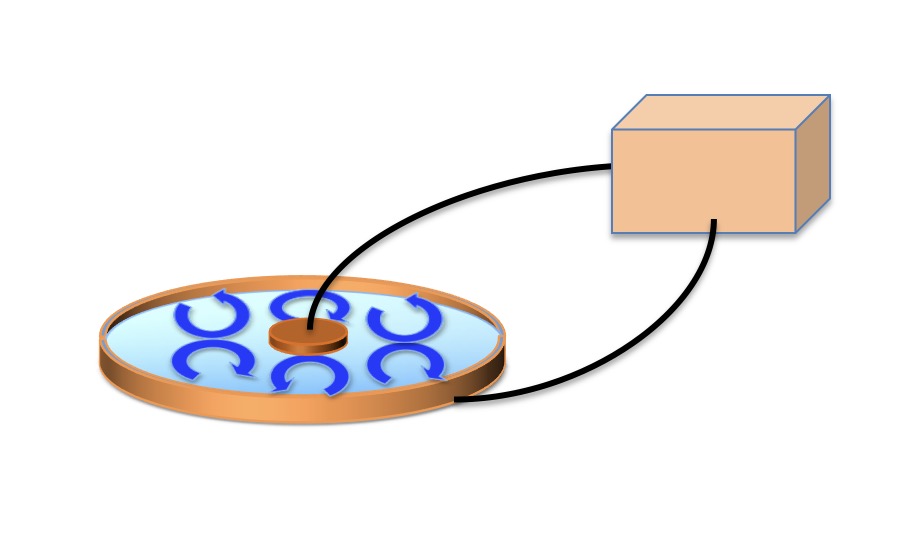} \\
				(a) & (b)
			\end{tabular}
		\end{center}
		\caption{Experimental settings. The liquid film is located in between two concentric electrodes. 
		(a) the two electrodes extend to infinity (infinite case) and
		(b) the two electrodes have negligible thickness (slim case). Notice that in both cases, the outer electrode is extended to infinity in the $xy$-plane (not pictured).}
		\label{fig:annular}
	\end{figure}

    In this work, we follow the experimental setting in \cite{daya1998annularshear} and \cite{TDDW2007electroconvection}.
    The thin film of liquid crystal is at room temperature and the crystals are arranged in layers with their long axis aligned with the normal of the film plane (semectic-A phase). 
    Under such arrangement, the fluid molecules cannot move across the film thickness, thereby \modif{restricted to} a planar motion.
    The fluid has a poor conductivity in the direction perpendicular to 
	the long axis of the liquid crystal molecules.
	As a consequence, the magnetic fields generated by the low current inside the fluid are neglected. 
	The thin-layered film is placed in between two concentric electrodes as illustrated in~\Cref{fig:annular}.
	Two different settings are considered: the electrodes extend to infinity in the \modif{fluid} normal direction (infinite case) and the electrodes have negligible thickness (slim case).

	The mathematical models are derived \modif{in the limiting case of vanishing thickness.} 
	The resulting surface charge density conservation relation and the mass and momentum conservation relations for the incompressible liquid are derived on a two dimensional bounded domain $\Omega$. 
	However, the electric field (or potential) remains \modif{described} in all $\mathbb R^3$. 
	To fully take advantage of a reduced modeling setting, we resort to equivalent nonlocal representations of the restriction to $\Omega$ of the electric potential.
	Depending on the electrodes configuration, these nonlocal representations involve either the \emph{spectral Laplacian} defined via spectral expansions (infinite electrodes) or the \emph{integral Laplacian} defined using a Fourier transform (slim electrodes).
	We refer to \cite{bonito2018numerical} for a survey of numerical methods associated with these two operators.

	Numerical simulations of electroconvection for slim electrodes are obtained in~\cite{TDDW2007electroconvection,Tsai2007thesis,jabbour2015simulation}.
	Second order backward differentiation schemes for the time evolution coupled with spectral methods for the space discretization of all the quantities  are advocated.  
	The space discretizations are based on standard continuous, piecewise polynomial finite elements: Taylor-Hood elements for the fluid dynamics and continuous piecewise (bi-)linear finite elements for the electric potential and the surface charge density. 
	At this point, it is worth mentioning that the discretization of the nonlocal operators only requires finite element solvers for standard reaction-diffusion elliptic partial differential equations, thereby making the overall algorithm easily implementable.
	The evolution of the velocity-pressure in the liquid is approximated using a second order backward difference method. 
	In our experimental settings and unlike for the fluid dynamics, the surface charge density conservation law \modif{may be convection dominated and prone to numerical instabilities.}
	\modif{To circumvent this issue, we propose} an explicit second order strongly stability preserving scheme (SSP) (see~\cite{gottlieb2003SSP}) coupled with a second order entropy viscosity stabilization from~\cite{Guermond2017invariant} for \modif{its  approximation.}
	
	We assess numerically the advantages of the two different electrodes configurations.
	We also provide a numerical study on the effect of three critical nondimensional parameters in the electroconvection process.
	They are the Rayleigh number $\mathcal R$ representing the ratio between the electric forcing and viscous dissipation, the Prandtl number $\mathcal P$ measuring the ratio between the charge relaxation time and the viscous relaxation time, and the circular electrodes aspect ratio $\alpha$.

	The outline of this paper is as follows. 
	We begin with the derivation of the full mathematical model in Section~\ref{s:mathematical} followed with its two dimensional model reduction in Section~\ref{s:model_reduction}.
	The numerical algorithms for the approximations of the fluid dynamics, the surface charge density and electric potential are described in Section~\ref{s:numerical}.
	Their general performances together with numerical assessments of  favorable conditions to obtain electroconvection are discussed in Section~\ref{s:simu}. 
	When using the same experimental configurations as in \cite{TDDW2007electroconvection,Tsai2007thesis}, our findings are in good agreement. This is also detailed  in Section~\ref{s:simu}.
	Conclusions are drawn in Section~\ref{s:conclusion}.

	\section{Mathematical models}\label{s:mathematical}
	
	We describe below the derivation of the electroconvection model from Maxwell system.
	We consider two settings: infinite and slim electrodes. 
	The former corresponds to the simulation reported in \cite{daya1998annularshear, daya1999annularshear, TDDW2007electroconvection} while the later is related to the analysis in \cite{CEIV2016ec}.
	Worth mentioning, vanishing charge densities on $\partial \Omega$ are considered in \cite{CEIV2016ec} allowing for smoother charge densities and simplifications of the mathematical model. 
	\revone{We do not make such assumption incompatible with the requirement of charge conservation and refer to Remark~\ref{r:constantine} for further details.}
	
	We recall that our experimental setting consists of a thin layer of liquid crystal in an annular region.
	\revone{We assume that the liquid crystal molecules can only move along the $xy-$plane and polarized in the $z$ direction without being able to change polarization. This corresponds to the so-called Smectic-A phase, which manifests itself when the liquid crystal are at room temperature $24\pm2^\circ$C.}
	\revone{At this point, we remark that the assumptions made on the state of the liquid crystals are critical for the subsequently derived model reduction. For instance, if instead in a Smectic-C phase,  the long axis of liquid crystal molecules pertain a fixed angle with $z$ direction.  The possibility of model reduction in this context is open for further investigations.}
	
	\subsection{Geometry}
	We denote by
	$\Omega:=\{\bx \in \RR^2\ :\  R_i < \abs{\bx} < R_o\}$, $0<R_i<R_o<\infty$,
	the annular region. 
	The liquid crystal is confined in the domain $\Omega_s := \Omega\times (-s,s)$, where $s>0$ stands the film thickness.
	We also denote by $K_i$ and $K_o$ the regions in $\Omega$ occupied by the inner and outer electrodes respectively:
	$K_i:= \{ \bx \in \RR^2 \ :  \ | \bx| < R_i \}$ and 
	$K_o:= \{ \bx \in \RR^2 \ :  \  | \bx| > R_o \}$.
	
	As already mentioned, we consider two experimental settings: 
	$K_i \times \{ 0 \}$ and $K_o \times \{ 0 \}$ (slim electrodes) 
	or $K_i^\infty:= K_i \times \mathbb R$ and $K_o^\infty:= K_o \times \mathbb R$ (infinite electrodes).
	Generically, we use the notation $\mathcal K_o$ and $\mathcal K_i$ to denote either the slim or infinite electrodes. 
	Furthermore, we \modif{set} $D_s:= \RR^3 \setminus (\overline{\Omega}_s\cup \overline{\mathcal K}_i \cup \overline{\mathcal K}_o)$, with $D:=D_0$, to denote the free space.
	
	\subsection{Electro-magnetism}
	The liquid crystals considered exhibit poor conductivity in the xy-plane, thereby resulting in low currents inside the liquid film and negligible magnetic field.
	Without magnetic effects, the electric field satisifes
	\begin{equation*}
	\nabla \times \bE = 0, 
	\end{equation*}
	which guarantees the existence of a potential function $\Psi: \mathbb R^3  \rightarrow \mathbb R$ satisfying
	\begin{equation*}
	\bE = -\nabla \Psi.
	\end{equation*}

	The relation between the potential $\Psi$ and  the surface charge density $q:\partial \Omega_s \rightarrow \mathbb R$ (units: \si{\coulomb\per\metre\tothe{2}}) is derived from the Gauss law. 
	Because there is no charge in the free space $D_s$, it directly implies that
	\begin{equation} \label{e:psi_free_full}
	\DIV \bE = 0 \qquad \textrm{or} \qquad \Delta \Psi = 0 \qquad \textrm{in} \quad D_s.
	\end{equation}
	In addition, the electric potential $\Psi$ is set to the appropriate voltages on the electrodes 
	\begin{equation}\label{e:psi_bc}
	\Psi = V \quad \textrm{on} \quad \mathcal K_{i} \qquad \textrm{and} \qquad \Psi = 0 \quad \textrm{on} \quad \mathcal K_{o}.
	\end{equation}
	\modif{Behavior of the electric potential at $z\to \infty$ is required to close the system.
	In the case of slim electrodes, one assumes that
	\begin{equation}\label{e:psi_decay_slim}
	\lim_{z \to \infty} | \Psi(x,y,z)| =0, \qquad (x,y) \in \mathbb R^2,
	\end{equation}
	while
	\begin{equation}\label{e:psi_decay_infty}
	\lim_{z \to \infty} | \Psi(x,y,z)| < \infty, \qquad (x,y) \in \Omega,
	\end{equation}
	is imposed in the infinite electrodes case. 
	}
	
\begin{figure} [ht]
		\centering
		\includegraphics[width=0.8\textwidth]{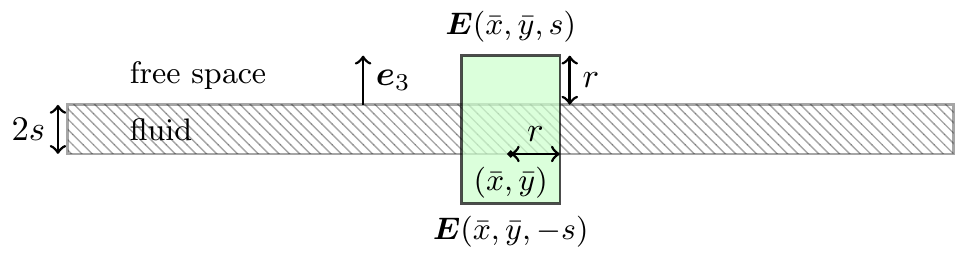}
		\caption{Gaussian pillbox enclosing a small region of the fluid. }
		\label{fig:gauss_box}
	\end{figure}

	For the liquid region $\Omega_s$, we recall that even with low conductivity, the charges are located at the boundary of $\Omega_s$. 	
	\modif{To derive a relation between the electric field and the surface charge density, we} consider a ``Gaussian pillboxes" (see~\Cref{fig:gauss_box} and also  \cite{griffiths2005electro})
	$$
	G(\bar x,\bar y,r):=(\bar x - r, \bar x +r) \times (\bar y - r,\bar y+r) \times (-s-r,s+r)
	$$
	centered at $(\bar x,\bar y) \in \Omega$ with $r =s$.
	Applying the Gauss law on $G(\bar x,\bar y,r)$ yields
	\begin{equation}\label{e:guass_pillbox}
	\varepsilon_0(\bE(\bar x,\bar y,s)  -\bE(\bar x,\bar y,-s))\cdot \be_3  = 2q + O(s).
	\end{equation}
	Here $\be_3 = (0,0,1)$ and $\varepsilon_0$ is the permittivity of the free space.	
	Notice that the constant 2 is due to contribution from the top and bottom sides of the film.
	Also, the permittivity of the liquid crystal is assumed to be isotropic in the $xy$ directions 
	and therefore the contributions from the other sides are of $O(s)$ in \eqref{e:guass_pillbox}.
	
	\subsection{Fluid Dynamics}\label{s:fluid_dyn}
	\revone{In the Smectic-A phase, the liquid crystal behaves like a two-dimensional incompressible fluid. The Navier-Stokes system is used to model the relation between the pressure $p(t):\Omega_s \rightarrow \mathbb R$ and the fluid velocity $\bu(t) :\Omega_s \rightarrow \mathbb R^3$}
	\begin{equation}\label{e:navierstokes}
	\rho ( \pt + \udg) \bu + \eta \Delta \bu + \nabla p 
	= \bL,  
	\quad \nabla \cdot \bu = 0, 
	\quad \textrm{both in }\Omega_s \times \RR_+,
	\end{equation}
	where $\rho$ is the fluid density, $\eta$ is the shear viscosity and $\bL := q \nabla \Psi$ is the Lorentz force induced by the motion of the electric current.
	The velocity \revone{is two dimensional $\bu \cdot \be_3 =0$ in $\Omega_s$ (and in particular on the top and bottom sides). 
	At the electrodes $\partial \Omega_s \cap \partial \mathcal K_{i/o}$, we restrict further the motion and consider a no-slip conditions $\bu = 0$.}
	
	\subsection{Small Thickness Limiting Model}
	We now consider the limiting model when the thickness $s$ tends to $0$.
	Relation~\eqref{e:guass_pillbox} between the electric field $\bE = - \nabla \Psi$ and the surface charge density $q$ becomes
	\begin{equation} \label{e:jump}
	\lim_{z\downarrow 0} \frac{\partial}{\partial z} \Psi - \lim_{z\uparrow 0} \frac{\partial}{\partial z} \Psi = -\frac{2q}{\varepsilon_0}
	\quad \textrm{in } \Omega,
	\end{equation}
	Notice that by symmetry, relations~\eqref{e:psi_free_full}, \eqref{e:jump} and \eqref{e:psi_bc} reduces to a system of partial differential equations on the half plane $D^+:= D \cap \{ z \geq 0 \}$.
	It reads 
	\begin{equation} \label{e:psi_free}
	\Delta \Psi = 0, \quad \In D^+,
	\end{equation}
	together with the boundary conditions
	\begin{equation}\label{e:psi_bc_all}
	\begin{split}
	& \Psi = V ~ \textrm{on} ~ \mathcal K_{i}^+, ~ \quad \Psi = 0 ~ \textrm{on} ~ \mathcal K_{o}^+, ~\textrm{and } \lim_{z \downarrow 0} \partial_z \Psi = -\frac{q}{\varepsilon_0} ~ \textrm{on} ~\Omega,
	\end{split} 
	\end{equation}	
	where $\mathcal K_{i/o}^+:= \mathcal K_{i/o} \cap \{ z \geq 0 \}$.
	
	\revone{Thanks to the two dimensional assumption on the fluid motion (Smectic-A phase), the Navier-Stokes system \eqref{e:navierstokes} reduces to}
	\begin{equation}\label{e:NS}		
\begin{split}	
\rho_\Omega ( \pt + \bu_\Omega \cdot \nabla_\Omega ) \bu_\Omega + \eta_\Omega \Delta_\Omega \bu_\Omega + \nabla_\Omega p_\Omega 
	& = \bL_\Omega,  \qquad \textrm{in }\Omega \times \RR_+,
\\
	\nabla_\Omega \cdot \bu_\Omega & = 0, \quad \qquad \textrm{in }\Omega \times \RR_+,
	\end{split}
	\end{equation}
	where $p_\Omega(t) := p(t)|_\Omega:\Omega  \rightarrow \mathbb R$, $\bu_\Omega(t): \Omega \rightarrow \mathbb R^2$ denotes the first two components of $\bu$ restricted to $\Omega$ and $\rho_\Omega$, $\eta_\Omega$ are the two-dimensional fluid mass density (units: \si{\kilogram\per\metre\tothe{2}}) and
	two-dimensional  fluid shear viscosity (units: \si{\pascal\second\meter}) respectively. 
	Also $\bL_\Omega$ is the projection of the Lorenz force $\bL$ on the plane supporting $\Omega$, i.e.  $\bL_\Omega = q \nabla_\Omega \Psi_\Omega$. 
	In addition, the boundary conditions on the velocity become
	\begin{equation}\label{e:NS_bc}
	\bu_\Omega = 0 \qquad \textrm{on}\quad \partial \Omega \times \mathbb R_+.
	\end{equation}

	In $\Omega$, the charge density flux and electric field are assumed to satisfy the Ohm's law, i.e., 
	across any closed curves with outside pointing normal $\bnu$ we have
	$$
	\nabla_\Omega q \cdot \bnu  = -  \sigma_\Omega  \bE_\Omega \cdot \bnu,
	$$
	where $\sigma_\Omega$ stands for the two-dimensional fluid electrical conductivity (units: \si{\siemens}) and $\bE_\Omega$ denotes the first two component of $\bE$ (projection onto the plan supporting $\Omega$).
	With this assumption, the conservation of surface charge density reads
	\begin{equation*}
	\pt q + \nabla_\Omega \cdot (\bu_\Omega q -  \sigma_\Omega \bE_\Omega) = 0, \qquad \textrm{in } \Omega \times \mathbb R_+.
	\end{equation*}
	or
	\begin{equation}
	\pt q + \nabla_\Omega \cdot (\bu_\Omega q)  -  \sigma_\Omega \Delta_\Omega \Psi_\Omega = 0, 
	\qquad \textrm{in } \Omega \times \mathbb R_+,
	\label{e:charge}
	\end{equation}
	using the electric potential $\Psi_\Omega := \Psi|_\Omega$.
	The surface charge density $q$ in \eqref{e:charge} is defined up to a constant fixed by the total charge conservation relation
	\begin{equation}\label{e:charge_mean}
	\int_\Omega q(t) = \int_\Omega q(0), \qquad t > 0.
	\end{equation}
	
	In summary, the electric potential $\Psi$, the surface charge density $q$, the velocity and pressure $(\bu,p)$ are related by the system of differential equations \eqref{e:psi_free}, \eqref{e:NS},  \eqref{e:charge}, and \eqref{e:charge_mean}, which is supplemented with the boundary conditions \eqref{e:psi_bc_all} and \eqref{e:NS_bc}.
	From now on, we drop the subscript $\Omega$ on $\bu_\Omega$, $p_\Omega$.
	
	\subsection{Nondimensional Model}
	To sort out the effects of the different parameters present in the model, we follow \cite{Tsai2007thesis} and rewrite the governing equations~\eqref{e:psi_free}, \eqref{e:jump}, \eqref{e:charge}, and \eqref{e:NS} using the rescaled variables 
	$$ 
	\hat{\bfx}  := \frac{\bfx}{d},\;
	\hat t := \frac{\sigma_\Omega}{\varepsilon_0 d} t,  \quad \widehat{\Omega} := \{ \hat{\bfx} \ :  \bfx  \in \Omega\}, \quad \textrm{and} \quad \widehat{D}^+ := \{  \hat{\bfx} \ :  \bfx  \in  D^+  \},
	$$
	where $d := R_o-R_i$ is the distance between the two electrodes.
	Similarly for the electrodes, we set $\widehat{\mathcal K}_{i/o}^+:=  \{  \hat{\bfx} \ :  \bfx  \in  \mathcal K_{i/o}^+  \}$.
	In addition, we define the rescaled functions
	\begin{align*}
	\widehat{\Psi}( \hat{\bfx}, \hat t) &:= \frac {1}{V}\Psi( \bfx,t), \quad	\widehat{\Psi}_{\widehat{\Omega}}=\widehat{\Psi}|_{\widehat{\Omega}},
	\qquad 
	\hat{q}( \hat{\bfx}, \hat t) := \frac{d}{\epsilon_0V} q(\bfx,t),\\
	\widehat{\bu} ( \hat{\bfx}, \hat t)&:= \frac{\epsilon_0}{\sigma_\Omega}\bu(\bfx,t),\qquad
	\hat p ( \hat{\bfx}, \hat t) := \frac{\epsilon_0^2}{\sigma_\Omega^2\rho_\Omega}p(\bfx,t).
	\end{align*}
	In order to simplify the notations, we omit the notation $\hat{.}$ and for now on only consider the rescaled variables.
	In particular, we write
	\begin{equation}\label{eq:dl-EC-charge}
	\pt {q} + \bu \cdot \nabla_\Omega  q  - \Delta_\Omega \Psi_\Omega  = 0, \qquad  
	\In {\Omega} \times \RR_+,
	\end{equation}
	with $\int_\Omega q(t) = \int_\Omega q(0)$ for the conservation of surface charge density relation,
	\begin{equation}\label{e:dl-jump}
	-\Delta \Psi = 0, 
	\In D^+ \times \RR_+, 
	\end{equation}
	with boundary conditions 
	\begin{equation}\label{e:dl-psi_bc_all}
	\begin{split}
	& \Psi = 1 ~ \textrm{on} ~ \mathcal K_{i}^+, ~ \Psi = 0 ~ \textrm{on} ~ \mathcal K_{o}^+, ~\textrm{and}~ \lim_{z \downarrow 0} \partial_z \Psi = -q ~ \textrm{on} ~\Omega
	\end{split} 
	\end{equation}
	\modif{(together with \eqref{e:psi_decay_slim} in the case of slim electrodes and \eqref{e:psi_decay_infty} in the case of infinite electrodes)}
	for the relations between the electric potential and the surface charge density, 
	and
	\begin{equation}\label{eq:dl-EC-momentum}
	\pt \bu + (\bu\cdot\nabla_\Omega) \bu - \calP\Delta_\Omega \bu + \nabla_\Omega p  
	= -\calR\calP q\nabla_\Omega\Psi_\Omega, \qquad 	\nabla_\Omega\cdot\bu = 0,
	\end{equation}
	in $\Omega \times \RR_+$
	with boundary conditions $\bu = 0$ on $\partial \Omega$ for the Navier-Stokes system.
	The two dimensionless parameters $\calP$ and $\cal R$ appearing in \eqref{eq:dl-EC-momentum} are the Prandtl and Rayleigh numbers.
	They are given by 
	\begin{align} \label{PandR}
	\calP := \frac{\epsilon_0\eta_\Omega}{\rho_\Omega\sigma_\Omega d}, \quad \And \quad
	\calR := \frac{\epsilon_0^2V^2}{\eta_\Omega\sigma_\Omega}.
	\end{align}
	The Prandtl number indicates the fluid viscous relaxation ability relative to its charge relaxation ability 
	while the Rayleigh relates  the electric with the dissipation forces. In Section~\ref{s:simu}, we propose a numerical study determining what range of parameters allows for electroconvection.
	We also introduce a geometric characteristic parameter $\alpha := R_i/R_o \in (0,1)$ so that
	\begin{equation} \label{e:ratio}
	R_i = \frac{\alpha}{1-\alpha}, \qquad \textrm{and} \qquad
	R_o = \frac{1}{1-\alpha}.
	\end{equation}
	As we shall see in~\Cref{s:simu}, this parameter affects the number of vortices during the electroconvection process.
	
	\section{Model Reduction}\label{s:model_reduction}
	
	The electroconvection system \eqref{eq:dl-EC-charge}-\eqref{eq:dl-EC-momentum} is mainly two dimensional (i.e. defined on $\Omega$) except for \eqref{e:dl-jump}, where the electric potential must be computed in the entire free space $D^+$. 
	However, notice that only its trace in $\Omega$ is required in  \eqref{eq:dl-EC-charge} and \eqref{eq:dl-EC-momentum}.
	This, together with the two different electrodes configurations (infinite and slim), is exploited in Sections~\ref{s:reduce:infinite} and~\ref{s:reduce:finite} to replace \eqref{e:dl-jump} with nonlocal problems in $\Omega$.

	\subsection{Infinite Electrodes} \label{s:reduce:infinite}
\begin{figure}[h]
		\begin{center}
			\includegraphics[width=0.45\textwidth]{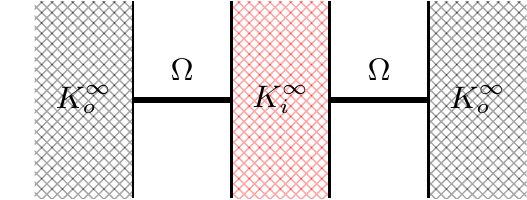} 
		\end{center}
		\caption{Domain and boundaries, cross section view in the infinite electrodes case; compare with Figure~\ref{f:slim_electrode}.}
		\label{fig:domain}
	\end{figure}

	We assume that the two electrodes $K_i$ and $K_o$ extend to infinity along the $z$ directions, refer to \Cref{fig:domain} for an illustration. 
	This is the setting considered for instance in \cite{CEIV2016ec}. However, because we do not impose vanishing charge densities on $\partial \Omega$, our model do not reduce to the one analyzed in \cite{CEIV2016ec}.
	Worse, it appears that less regular charge densities are to be expected (see Remark~\ref{r:constantine}).
	
	We decompose $\Psi$ in two parts $\Psi = \Psi_0 + \Psi_K$.
	The component $\Psi_K$ accounts for the voltage set by the electrodes and is defined as the solution to 
	\begin{equation} \label{eq:homo}
	\Delta \Psi_K = 0 	\In   D^+ \quad
	\end{equation}
	together with $\Psi_K  = 1$ on $\mathcal K_i^+$, $\Psi_K= 0$ on $\mathcal K_o^+$ and
	$ \lim_{z \downarrow 0} \partial_z \Psi_K = 0$ on $\Omega$. 
	Notice that $\Psi_K$ is independent of the $z$ variable.
	In fact, its exact expression is given by
	\begin{equation} \label{eq:bdd_psi_1}
	\Psi_K(x,y,z)  = \eta(x,y) : = \frac{ \ln (\sqrt{x^2+y^2}/(1-\alpha))}{\ln (\alpha)}.
	\end{equation}
	
	The second part $\Psi_0 = \Psi-\Psi_K$ depends on the charge density $q$ and solves
	\begin{equation*}
	-\Delta \Psi_0 = 0 \In D^+
	\end{equation*}
	together the boundary conditions  
	$$
	\Psi_0 = 0~\textrm{on}~ \mathcal K_i^+ \cup \mathcal K_o^+, \quad \textrm{and} \quad \lim_{z\downarrow 0} \partial_z \Psi_0 = -q \quad \textrm{on }\Omega.
	$$
	
	Following  \cite{ST2010extension}, we realize that the trace $\Psi_{0,\Omega}:= \Psi_0|_
	\Omega$ satisfies the following non-local partial differential equation on $\Omega$
	\begin{equation} \label{eq:bdd_psi_2}
	(-\Delta_\Omega)^{\frac 1 2} \Psi_{0,\Omega} = q \quad  \In \quad \Omega,
	\end{equation}
	where $(-\Delta_\Omega)^{\frac 1 2}$ is the \emph{spectral Laplacian} defined using a spectral expansion as \modif{described now}. We denote by $\{\phi_j\}_{j=1}^\infty \subset H^1_0(\Omega)$  
	an $L_2(\Omega)$-orthonormal basis of eigenfunctions of $-\Delta$ and by $\{\lambda_j\}_{j=1}^\infty$ the associated eigenvalues. 
	The spectral Laplacian is then defined for $v \in H^1_0(\Omega)$ via the relation
	\begin{equation}\label{spec_op}
	(-\Delta_\Omega)^{\frac 1 2} v := \suminf \lambda_j^\frac{1}{2} \left( \int_\Omega v \phi_j \right) \phi_j.
	\end{equation}
	
	Returning to the expression for $\Psi$, we find that its restriction to $\Omega$, $\Psi_\Omega:=\Psi|_\Omega$ satisfies
	\begin{equation}\label{e:Psi_trace_bounded}
	\Psi_\Omega = \Psi_{0,\Omega} + \eta
	\end{equation}
	where the expression of $\eta$ is given in ~\eqref{eq:bdd_psi_1}.
	The reduced model system, defined on $\Omega$, consists of \eqref{eq:dl-EC-charge}, \eqref{eq:dl-EC-momentum}, \eqref{eq:bdd_psi_2} and \eqref{e:Psi_trace_bounded}.
	
	\begin{remark}\label{r:constantine}
		We have already mentioned that compared with \cite{CEIV2016ec}, we do not impose $q|_{\partial \Omega} =0$ but rather
		$\int_\Omega q = 0$. 
		In particular, the charge density does not belong to $H^1_0(\Omega)$, the domain of the operator $(-\Delta_\Omega)^{\frac 12}$, and we cannot use the  relation  
		$$
		-(\nabla_\Omega \cdot \nabla_\Omega) (-\Delta_\Omega)^{-\frac 1 2}   = (-\Delta_\Omega) (-\Delta_\Omega)^{-\frac 1 2} =  (-\Delta_\Omega)^{\frac 1 2},
		$$
		to simplify the charge relation~\eqref{eq:dl-EC-charge} as
		$$
		\pt {q} + {\bu} \cdot \nabla_\Omega  q + (-\Delta_\Omega)^{\frac 1 2} q  = 0  
		\quad \In \quad {\Omega} \times \RR_+.
		$$
		The former is the starting point of the analysis proposed in~\cite{CEIV2016ec}.
		We did not make such assumption because it implies 
		$$
		\lim_{t\to \infty} \int_\Omega q(t) = 0,
		$$
		which is incompatible with the surface charge density conservation required in our context.
	\end{remark}
	\begin{remark}
		The decomposition \eqref{e:Psi_trace_bounded} of $\Psi_\Omega$ happens to correspond to the definition of the fractional laplacian with non vanishing boundary condition proposed in \cite{antil2017fractional}.
	\end{remark}

	\subsection{Slim electrodes}\label{s:reduce:finite}
\begin{figure}[h]
		\begin{center}
			\includegraphics[width=0.7\textwidth]{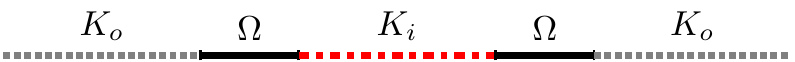} 
		\end{center}
		\caption{Domain and boundaries, cross section view in the slim electrodes case; compare with Figure~\ref{fig:domain}.}
		\label{f:slim_electrode}
	\end{figure}

	Instead of assuming infinite electrodes in the $z$ direction, we consider electrodes with negligible height as illustrated in~\Cref{f:slim_electrode}.
	
	In this setting as well, we separate the voltage and charge contributions in the electric potential 
	\begin{equation}\label{e:decomp_slim}
		\Psi = \Psi_0 + \Psi_K.
	\end{equation}
	In view of \eqref{e:dl-jump}, \eqref{e:dl-psi_bc_all} and the decay relation $\lim_{z\to +\infty} \Psi = 0$, the first part $\Psi_K$ is defined as the solution of
	\begin{equation} \label{e:harmo}
	\Delta \Psi_K = 0 \quad \In \RR^2\times \RR_+
	\end{equation}
	together with the boundary condition
	\begin{equation} \label{e:harmo_bc1}
	\Psi_K = \eta, \qquad \textrm{on }\mathbb R^2,
	\end{equation}
	and decay condition
	\begin{equation} \label{e:harmo_bc2}
		\lim_{z\to +\infty} \Psi_K = 0.
	\end{equation} 
	Here $\eta:\mathbb R^2 \rightarrow \mathbb R$ matches the voltage imposed by the electrodes and is extended harmonically in $\Omega$:
	$$
	-\Delta_\Omega \eta = 0 \quad \textrm{in }\Omega, \quad \eta = 0 \quad \textrm{on }K_o, \quad \textrm{and} \quad \eta = 1 \quad \textrm{on }K_i.
	$$
	Its exact expression on $\Omega$ matches \eqref{eq:bdd_psi_1}.
	
	The second component $\Psi_0 = \Psi-\Psi_K$ satisfies
	\begin{equation} \label{eq:nonhomo_unbdd}
	\begin{split}
	\left\{
	\begin{array}{rll}
	-\Delta \Psi_0 &= 0,  				&\In D^+,\\ 
	\partial_z\Psi_0 &= -q - \partial_z\Psi_K,   &\On \Omega,\\
	\Psi_0 &= 0, 						    &\On K_i \cup K_o,
	\end{array} 
	\right.
	\end{split}
	\end{equation}
	together with the decay condition $\lim_{z \to +\infty}\Psi_0 = 0$.
	It depends on the value of $ \partial_z\Psi_K|_\Omega$ we determine now.
	Standard techniques based on a separation of variables  $\Psi_K(x,y,z) =X(x,y) Z(z)$ and the representation of  $X(x,y)$ in terms of Bessel's functions of the first kind (\cf ~\cite[Chapter~4.1]{Evans2010pde} and \cite{Tsai2007thesis}) reveal that together with the axy-symmetry property of $\Psi_K|_\Omega = \eta$, the solution to~\eqref{e:harmo}, \eqref{e:harmo_bc1}, \eqref{e:harmo_bc2} expressed in cylindrical coordinates $(r,\theta,z)$ takes the form
	\begin{equation}\label{e:PsiK}
	\Psi_K (r,\theta,z) 
	=  \int_{0}^{\infty} 
	e^{-kz}J_0(kr) A_0(k)dk.
	\end{equation}
	The function $A_0$ is determined using the boundary condition $\Psi_K(r, \theta,0) = \eta(r)$ and the orthogonality properties of $J_0$:
	$$
	A_0(k) = k \int_0^\infty \eta(r) J_0(kr)r dr = \frac{1}{k^2 \ln \alpha} \left( J_0(R_0 k) - J_0(R_i k)\right),
	$$
	where we used the expression \eqref{eq:bdd_psi_1} of $\eta$ to derive the second equality.
	We now take advantage of Formula~6.574.1 in \cite{GR2014table} to compute
	\begin{equation*}
	\begin{split}
	\lim_{z \downarrow 0} \partial_z \Psi_K(r,\theta,z)
	&= - \int_{0}^{\infty} k J_0(kr) A_0(k)dk \\
	& =- \frac{1}{\ln \alpha} \int_0^{\infty}  
	(J_0(R_ok) - J_0(R_ik))J_0(rk)dk \\
	& =- \frac{1}{\ln \alpha} \left(\frac{1}{R_o}
	{}_2F_1(\frac{1}{2},\frac{1}{2};1;\frac{r^2}{R_o^2}) 
	- \frac{1}{r} {}_2F_1(\frac{1}{2},\frac{1}{2};1;\frac{R_i^2}{r^2}) \right) 
	\end{split}
	\end{equation*}
	where the function ${}_2F_1$ denotes the hyper-geometric function, see \cite[Chapter~15]{Abramowitz1964handbook}. 
	In short, we write
	$$
	\lim_{z \downarrow 0} \partial_z \Psi_K = g,
	$$
	where
	\begin{equation}\label{d:g}
	g(r):= - \frac{1}{\ln \alpha} \left(\frac{1}{R_o}
	{}_2F_1(\frac{1}{2},\frac{1}{2};1;\frac{r^2}{R_o^2}) 
	- \frac{1}{r} {}_2F_1(\frac{1}{2},\frac{1}{2};1;\frac{R_i^2}{r^2}) \right).
	\end{equation}
	Notice that $g \in L^2(\Omega)$ thanks to the property
	\begin{equation*}
	\lim_{\beta \to 1^-} \frac{{}_2F_1(a,b;a+b;\beta)}{-\ln(1-\beta)} 
	= \frac{\Gamma(a+b)}{\Gamma(a)\Gamma(b)}.
	\end{equation*}
	As a consequence, the system~\eqref{eq:nonhomo_unbdd} for $\Psi_0$ is uniquely determined by $q$ and $g$. 
	
	We denote by $\Psi_{0,\Omega}$, the trace of $\Psi_0$ on $\Omega$, and define $\widetilde \Psi_{0,\Omega}$ to be its zero extension to $\mathbb R^2$, i.e.
	$$
	\widetilde \Psi_{0,\Omega} = \Psi_0\qquad \textrm{on}\quad \Omega, \qquad \textrm{and} \qquad \widetilde \Psi_\Omega = 0 \quad \textrm{on} \quad \mathbb R^2 \setminus \Omega.
	$$
	It turns out that $\widetilde \Psi_{0,\Omega}$ satisfies the following nonlocal problem
	\begin{equation} \label{eq:Psi0}
	(-\Delta_{\mathscr{F}})^{\frac 1 2} \widetilde{\Psi}_{0,\Omega} = q + g \qquad \textrm{in}\quad \Omega, 
	\end{equation}
	where the fractional operator $(-\Delta_\mathscr{F})^{\frac 1 2}$ is the integral fractional laplacian defined via the Fourier transform $\mathscr{F}$ in $\mathbb R^2$ as detailed now.
	For $v \in H^1(\mathbb R^2)$, the integral fractional Laplacian is a
	pseudo-differential operator with symbol $\abs{\zeta}$, i.e.
	\begin{equation}\label{inte_op}
	\mathscr{F} ((-\Delta_{\mathscr{F}})^{\frac 1 2} v)(\zeta) := \abs{\zeta}\mathscr{F}(v)(\zeta).
	\end{equation}
	
	\revone{The claim that the restriction to $\Omega$ of the solution to~\eqref{eq:nonhomo_unbdd} is the unique solution to~\eqref{eq:Psi0} directly follows from an argument provided in \cite{CS2007extension}, where $\partial_z\Psi_0$ is imposed on $\mathbb R^2$ rather than only in $\Omega$ as in our context. 
	Following the argument provided in Section 3.2 in  \cite{CS2007extension}, it is easy to see that
	the unique solution $\Psi$ to  \eqref{eq:nonhomo_unbdd} minimize the energy
	$$
	\frac 1 2 \int_{\mathbb R^2} \int_0^\infty | \nabla \Phi (x,y,z)|^2 ~dz~dxdy -  \int_{\Omega} (q(x,y)+g(x,y)) \Psi(x,y,0) ~dxdy
	$$
	over all functions $\Phi:\mathbb R^3 \rightarrow \mathbb R$ vanishing in $\mathbb R^3 \setminus \Omega$ and when $z \to \infty$.
	Furthermore, its restriction  $\Phi(z,y,0)$ also minimizes the energy
	$$
	J(w):= \frac 1 2 \int_{\mathbb R^2} |\xi| | \mathscr{F}(w)|^2 d\xi  \ -  \frac 1 \pi \int_{\Omega} (q(x,y)+g(x,y)) w(x,y) ~dxdy
	$$
	over all functions $w:\mathbb R^2 \rightarrow \mathbb R$ vanishing on $\mathbb R^2 \setminus \Omega$. 
	Since the solution to \eqref{eq:Psi0} is the unique minimizer of $J(.)$, the claim follows. 
	Details are available in the Appendix A of \cite{peng}.
}
	
	Combining the relations obtained for $\Psi_K$ and $\Psi_0$, we find that the restriction of $\Psi_\Omega = \Psi|_\Omega$ is given by
	\begin{equation} \label{e:Psi_trace_unbdd}
	\Psi_\Omega = \Psi_{0,\Omega}  + \eta
	\end{equation}
	where $\Psi_{0,\Omega}$ satisfies~\eqref{eq:Psi0} with $g$ as in \eqref{d:g} and $\eta$ is given by \eqref{eq:bdd_psi_1}. 
	
	The reduced model for the slim electrodes case defined on $\Omega$ consists of \eqref{eq:dl-EC-charge}, \eqref{eq:dl-EC-momentum} and \eqref{e:Psi_trace_unbdd}.
	
	\section {Numerical Algorithms}\label{s:numerical}
	
	In this section we detail the numerical algorithms advocated to approximate the fluid dynamic~\eqref{eq:dl-EC-momentum}, the surface charge density convection~\eqref{eq:dl-EC-charge}, and the two non-local problems for the electric field~\eqref{e:Psi_trace_bounded} or~\eqref{e:Psi_trace_unbdd} depending on the assumption made on the electrodes.	
	In fact, the time marching algorithm proposed consists of three sub-steps.
	First, the electric potential is approximated in $\Omega$ using the (previous) surface charge density.
	Second, the surface charge density approximation is updated using the electric potential and (previous) fluid velocity.
	Third, the fluid velocity (and pressure) is  updated with a Lorentz force computed using the surface charge density and electric potential.
	We detail each step separately below.
	Notice that the initial surface charge density and velocity are given allowing the algorithm to start.

	\subsection{Approximation of the Electric Potential~\eqref{eq:bdd_psi_2} or~\eqref{eq:Psi0} }\label{subsec:pot}

	We discuss the two different configurations separately.
	 	
	\paragraph{Infinite Electrodes} We start with the simpler case of infinite electrodes, where the electric potential satisfies~\eqref{eq:bdd_psi_2} involving the spectral Laplacian.
	Recall that $\Psi_\Omega = \Psi_{0,\Omega} + \eta$ \modif{with $\eta$ given} by~\eqref{eq:bdd_psi_1}.
	
	We adopt the numerical procedure put forward in \cite{BP2015elliptic,BP2016accretive}, which relies on the integral representation (Balakrishnan formula) of the electric potential $\Psi_{0,\Omega}$
	\begin{equation} \label{Balakrishnan}
	\Psi_{0,\Omega} = \frac{2}{\pi}  \int_{-\infty}^{\infty} e^{s} \Phi(s;q) \,d s.
	\end{equation}
	Here $\Phi(s;q)$ solves
	\begin{equation}\label{e:Phi_sub}
	\Phi(s;q)-e^{2s}\Delta \Phi(s;q) 
	=  q \qquad \textrm{in }\Omega, \quad \textrm{with}\quad \Phi(s;q)=0 \qquad \textrm{on }\partial \Omega.
	\end{equation}
	
	An exponentially convergent sinc quadrature, see \cite{BP2015elliptic,BP2016accretive,BLP2017accretive2}, is used to approximate the undefined integral.
	It reads
	\begin{equation} \label{e:bdd_sinc}
	\Psi_{0,\Omega} \approx \Psi^{k} 
	:=  \frac{2}{\pi} k \sum_{j = -N_{-}}^{N_{+}} e^{s_j} \Phi(s_j;q),
	\end{equation}
	where $k>0$ is the quadrature step size, $s_j := jk$, 
	$N_{+} := \left\lceil \frac{\pi^2}{k^2} \right\rceil, 
	\And  
	N_{-} := \left\lceil \frac{\pi^2}{k^2} \right\rceil$. 
	In the numerical simulations proposed in Section~\ref{s:simu}, we fix $k = 0.04$ so that $N_{-} = N_{+} = 62$.

\begin{figure}[ht]
		\begin{center}
			\begin{tabular}{>{\centering\arraybackslash} m{5cm} >{\centering\arraybackslash} m{5cm}}
				\includegraphics[scale=.08]{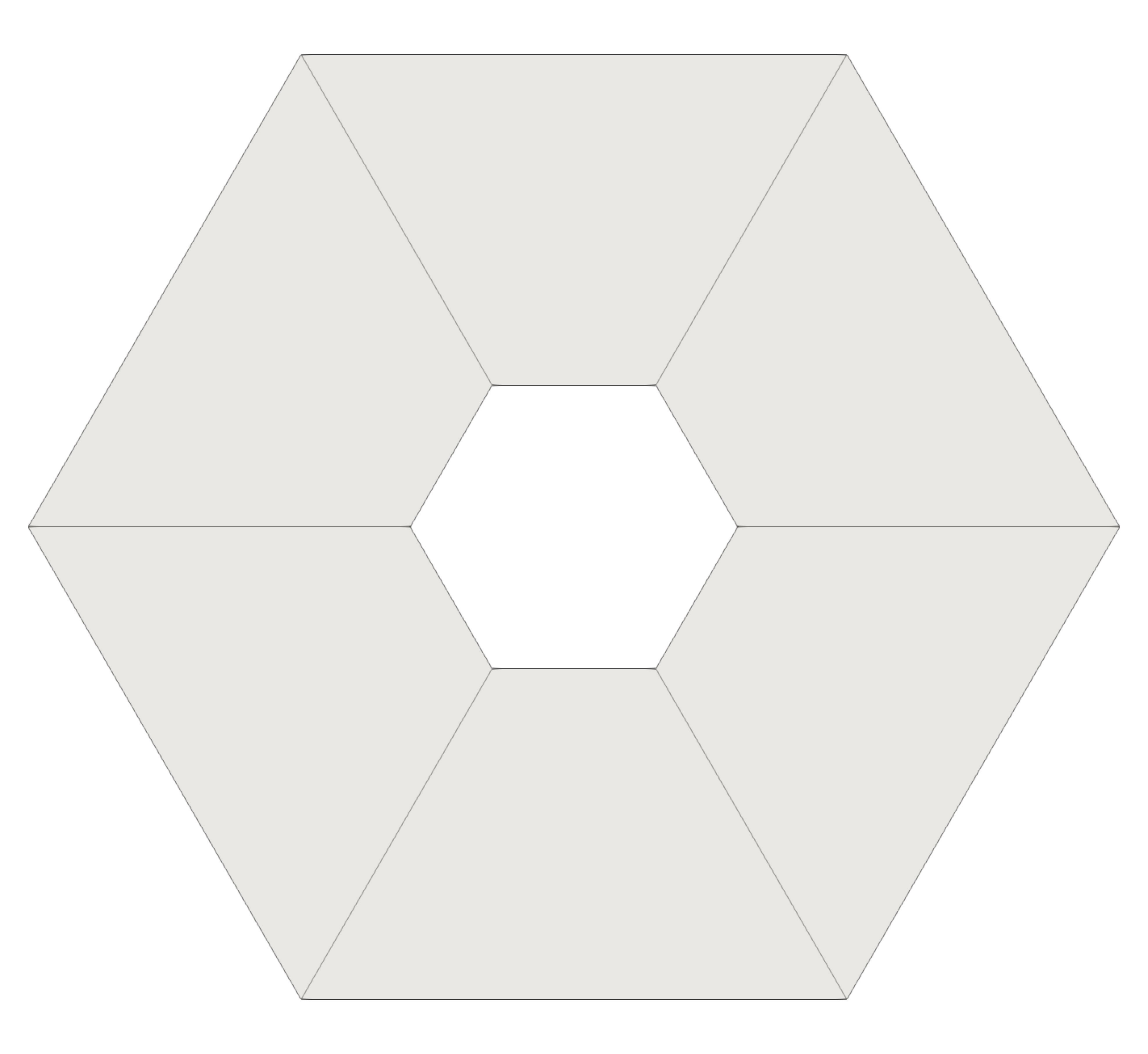}&
				\includegraphics[scale=.09]{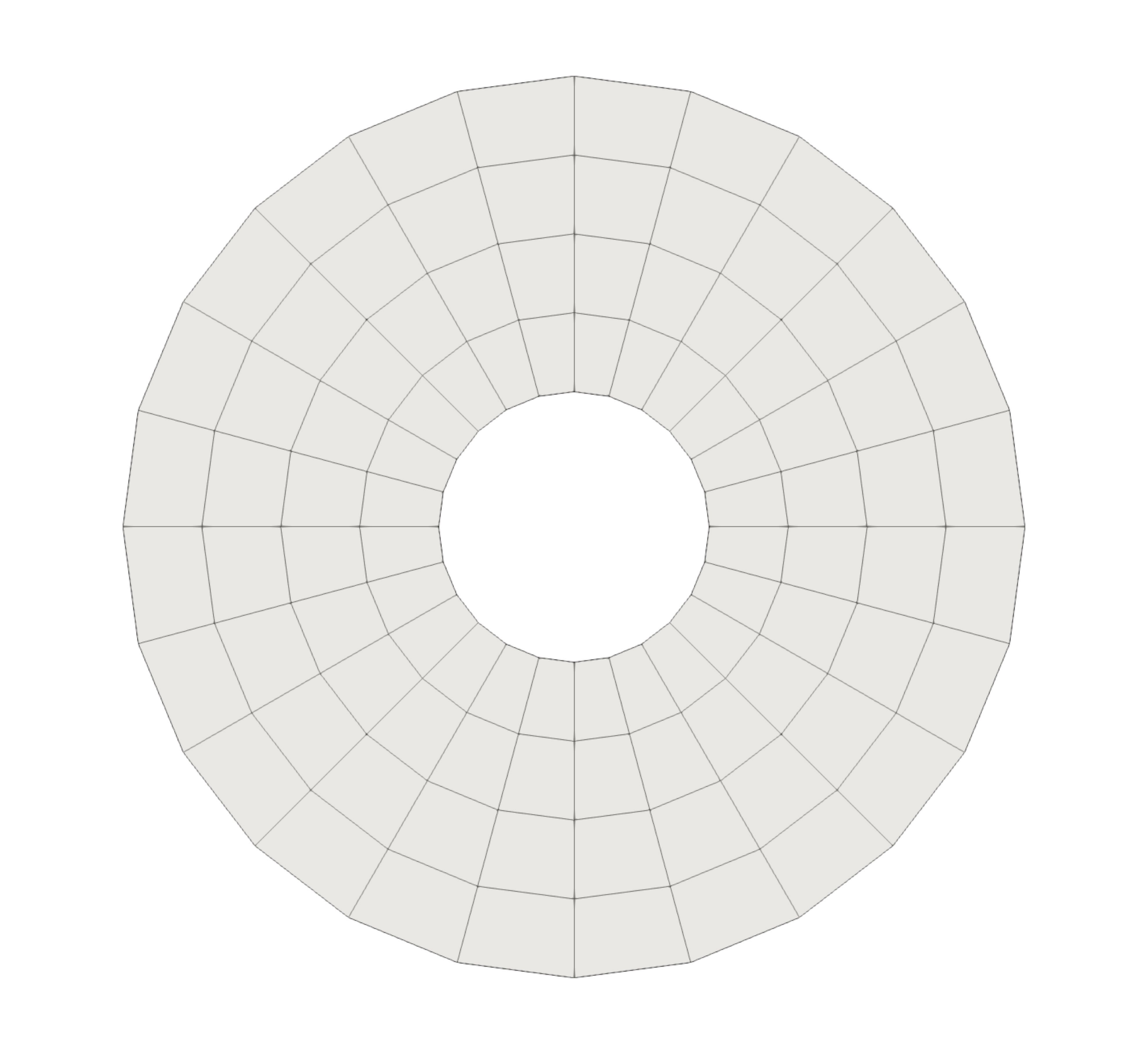} \\
				(a) & (b)
			\end{tabular}
		\end{center}
		\caption{Polygonal approximation $\Omega_h$ of $\Omega$ for an aspect ratio $\alpha = 0.33$. 
		(a) The coarse initial subdivision;
		(b) Approximation resulting in two successive uniform refinements and placing the boundary vertices on the boundary of $\Omega$.}
		\label{fig:infinite_mesh}
	\end{figure}

	The solutions $\Phi(s_j,q)$, $j=-N_-,...,N_+$ of the subproblems~\eqref{e:Phi_sub} are in turn approximated by continuous piecewise linear functions. We start with a coarse polygonal approximation of $\Omega$  subdivided into quadrilaterals as in Figure~\ref{fig:infinite_mesh}(a). 
	This polygonal approximation is then uniformly refined using quad-refinement but placing the boundary vertices on the exact boundary of $\Omega$. This gives rise to a polygonal domain $\Omega_h$ and a partition $\mathcal T_h$ made of quadrilaterals without hanging nodes, see for instance Figure~\ref{fig:infinite_mesh}(b). 
	Here $h$ denotes the maximum diameter of elements on $\mathcal T_h$.  
	
	The number of quadrilaterals in the coarse subdivision depends on the aspect ratio in order to maximize the quality of the subdivision, see \url{GridGenerator::hyper_shell} documentation in \cite{dealII85}.
	Furthermore, the final resolution $h$ required for the simulations presented in Section~\ref{s:simu} depends on the aspect ratio $\alpha$ as well. 
	Indeed, we shall see in Section~\ref{subsec:geometry} that larger aspect ratios yield more pairs of vortices during electroconvection, thereby increasing the resolution required.  The number of uniform refinements performed on the coarse subdivisions for each aspect ratio $\alpha$ are given in Table~\ref{t:infinite_coarse}.
	\begin{table}[ht!]
\begin{center}
		\begin{tabular}{c|ccccccc}
			$\alpha$ & 0.1	& 0.2  & 0.33 & 0.452 & 0.56 & 0.6446 & 0.9 \\
\hline
			$\# \mathcal T$ & 4 & 5 & 6 & 9 & 12 & 15 & 29 \\
			$\#$ uniform ref. & 5 & 5 & 5 & 5 & 5 & 5 & 6
			\end{tabular}
\end{center}	
	\caption{Infinite electrodes configuration: Number of quadrilaterals $\# \mathcal T$ and subsequent uniform refinements used in the coarse subdivision for different values of aspect ratio $\alpha$.}
	\label{t:infinite_coarse}
	\end{table}

	For each $T \in \mathcal T_h$, we denote by $F_T: \lbrack 0,1\rbrack^2 \rightarrow T$ the bilinear (invertible) map from the unit square to $T$. This map is instrumental to define the finite element space
	$$
	V_h:= \{ v \in C^0(\overline{\Omega}_h)\ : \ v|_T \circ F_T \textrm{ is bilinear for all } T \in \mathcal T_h \quad \textrm{and} \quad v|_{\partial \Omega_h}=0 \}.
	$$
	With these notations, the finite element approximation of $\Psi_{0,\Omega}^k$ in \eqref{e:bdd_sinc} is given by
	\begin{align} \label{eq:fem_direct}
	\Psi_{0,h} :=
	\Psi^{k}_{0,h}  
	:= \frac{2k}{\pi} \sum_{j = -N_{-}}^{N_{+}} e^{s_j} \Phi_h(s_j;q) \in V_h,
	\end{align}
	where  $\Phi_h(s_j;q) \in V_h$ solves
	\begin{equation}\label{e:Phi_h}
	\int_\Omega \Phi_h(s_j;q)  \varphi_h + e^{2s_j}\int_\Omega \nabla \Phi_h(s_j;q) \cdot \nabla \varphi_h	= \int_\Omega q  \varphi_h, 
	\Forall \varphi_h\in V_h.
	\end{equation}
	
	Returning to the decomposition~\eqref{e:Psi_trace_bounded}, we arrive at  
	\begin{equation}\label{e:final_Psi_infinite}
	\Psi_h:= \Psi_{0,h} + \pi_h \eta \in V_h,
	\end{equation}
	where $\pi_h$ stands for the $L^2$ projection onto $V_h$.
	
	It is worth mentioning, that the $N_{-}+N_{+}+1$ finite element problems \eqref{e:Phi_h} are mutually \textit{independent} making the parallel implementation straightforward. 
	Also, the algorithm consists of an outer loop \eqref{eq:fem_direct} gathering the contributions of the finite element solutions at each quadrature point $s_j$. 
	In particular, it only requires the implementation of a classical finite element solver for diffusion-reaction problems. 
	
	\paragraph{Slim electrodes} We now consider the case of slim electrodes.
	In this case, the electric potential on $\Omega$ is given by $\Psi_\Omega = \Psi_{0,\Omega} + \eta$ where $\Psi_{0,\Omega}$ satisfies~\eqref{eq:Psi0}.
	Unfortunately, an integral representation of $\Psi_{0,\Omega}$ like~\eqref{Balakrishnan} is not available for the integral Laplacian.
	Instead, we follow \cite{BLP2017integral} and take advantage of an integral representation of the action of $(-\Delta_{{\mathscr{F}}})^{\frac{1}2}$ by multiplying \eqref{eq:Psi0} by a smooth function compactly supported in $\Omega$ and integrating over $\Omega$. Theorem~4.1 in~\cite{BLP2017integral} guarantees that $\Psi_{0,\Omega}$ satisfies
\begin{equation}\label{e:phi_slim}
	\frac{1}{\pi} \int_{-\infty}^{\infty} e^{\frac{s}{2}} 
	\left(\int_\Omega (\Phi(s;\Psi_{0,\Omega})+\Psi_{0,\Omega}) \varphi\right) ds = \int_\Omega (q+g)\varphi
\end{equation}
for all smooth functions $\varphi$ compactly supported in $\Omega$. Here
$\Phi(s;\Psi_{0,\Omega}) \in H^1(\mathbb R^2)$ solves for all  $\varphi \in H^1(\RR^2)$
	\begin{equation}\label{e:full_eqs}
	e^s \int_{\mathbb R^2} \Phi(s;\Psi_{0,\Omega})  \varphi + \int_{\mathbb R^2} \nabla \Phi(s;\Psi_{0,\Omega}) \cdot \nabla\varphi  
	= - e^s\int_{\Omega} \Psi_{0,\Omega} \varphi.
	\end{equation}
	
	Similarly to the infinite electrode case, we can now use a sinc quadrature and standard finite element methods to provide an approximation of~\eqref{e:phi_slim}.
	The additional caveat stem from the fact that the domain of integration in \eqref{e:full_eqs} is $\mathbb R^2$.
	To circumvent this issue, we recall that 
	$$
	\Omega = \{ (x,y) \in \mathbb R^2 \ : \    \frac{\alpha}{1-\alpha} < \sqrt{x^2+y^2} < \frac{1}{1-\alpha} \}
	$$
	 and truncate the domain of integration to the disks
$$
	\Omega^M(s) := 
	\left\{
	\begin{array}{ll}
	\{(1+e^{-\frac{s}{2}}(M+1))(x,y)\;:\;  \sqrt{x^2+y^2} <1/(1-\alpha) \}, & e^{-\frac{s}{2}} \ge 1,\\ 
	\{(M+2)(x,y) \;:\;  \sqrt{x^2+y^2} <1/(1-\alpha)\}, & e^{-\frac{s}2} < 1,
	\end{array} 
	\right.
$$
	where $M>1$ is a given truncation parameter. 
	
	We postpone for the moment the discussion regarding the automatic subdivision of $\Omega^M(s)$ but denote by $\Omega^M_h(s)$ its polygonal approximation and by $V^M_h(s)$ the associated finite element space based on continuous piecewise bilinear finite elements vanishing on $\partial \Omega^M_h(s)$. 
	The approximation $\Psi_{0,h} \in V_h$ of $\Psi_{0,\Omega}$ satisfying \eqref{e:phi_slim} is then given by the relations
	$$
	\frac{k}{\pi} \sum_{j = -N_{-}}^{N_{+}} e^{\frac{s_j}{2}} 
	\int_{\Omega_h} (\Phi_{h}(s_j;\Psi_{0,h})+\Psi_{0,h}) \varphi_h = \int_{\Omega_h} (q+g)\varphi_h,\qquad \forall \varphi_h \in V_h,
$$
	where $\Phi_h(s_j;\Psi_{0,h}):=\Phi_{h}^{k,M}(s_j;\Psi_{0,h}) \in V_h^M(s_j)$ solves
	\begin{equation*}\label{e:full_eqs-trunc}
	\begin{split}
	e^{s_j} \int_{\Omega_h^M(s_j)} \Phi_{h}(s_j;\Psi_{0,h}) \varphi_h 
	+ \int_{\Omega^M_h(s_j)}\nabla \Phi_{h}(s_j;\Psi_{0,h}) &\cdot \nabla\varphi_h\\
	  &= - e^{s_j} \int_{\Omega_h} \Psi_{0,h}\varphi_h,
	\end{split}
	\end{equation*}
	for all $\varphi_h \in V_h^M(s_j)$.
	
	At this point, a few comments are in order. 
	The numerical procedure described above yields an approximation of the bilinear form~\eqref{e:phi_slim} associated with the relation~\eqref{eq:Psi0} for the electric potential but not the approximation $\Psi_{0,h}$ of $\Psi_{0,\Omega}$ directly.
	However, this is sufficient in the context of residual based iterative methods (such as conjugate gradient implemented in the present context), where only the action of the operator is required to compute residuals.
	The accumulation point of the iterative algorithm is $\Psi_{0,h} \in V_h$.
	
	For the numerical experiments provided in Section~\ref{s:simu}, we take $k=0.04$, $N_-=62$, $N_+=124$ and $M=3$.
	 It is shown in \cite{BLP2017integral} the proposed numerical method to approximate $\Psi_{0,\Omega}$ is exponentially convergent in the sinc quadrature parameter $k$ and in truncation parameter $M$.
	It is of optimal order (depending on the regularity of the exact electric potential $\Psi_{0,\Omega}$) in the space discretization.
	 We refer to \cite{BLP2017integral} for additional information.
	
	Returning to the decomposition \eqref{e:Psi_trace_unbdd}, we define
	\begin{equation}\label{e:final_Psi_slim}
		\Psi_h:= \Psi_{0,h} + \pi_h \eta \in V_h,
	\end{equation}
	where $\eta$ is given by~\eqref{eq:bdd_psi_1}.
	
\begin{figure}[h]
		\begin{center}
			\begin{tabular}{cc}
				 \includegraphics[scale=.09]{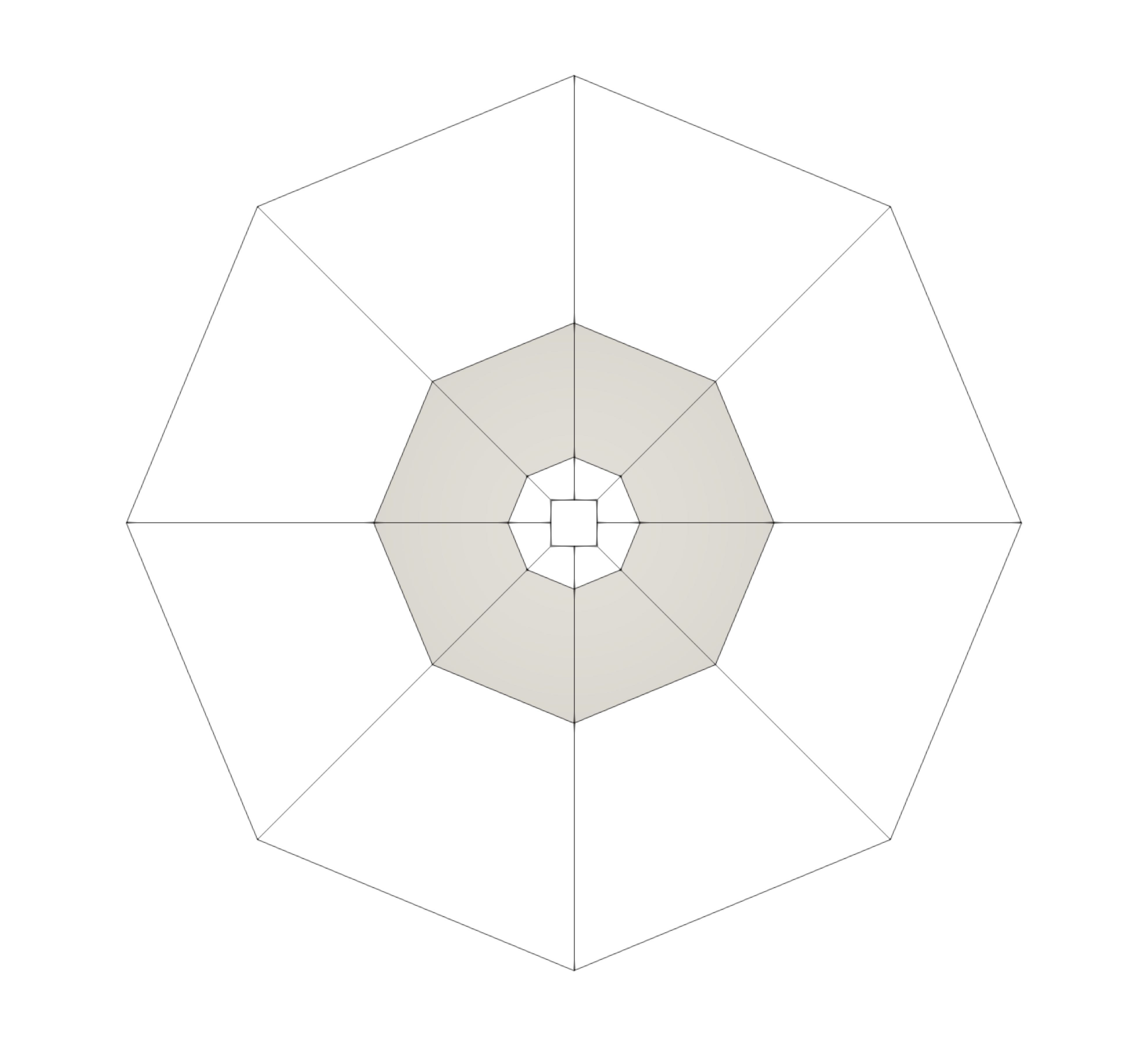} &
				 \includegraphics[scale=.09]{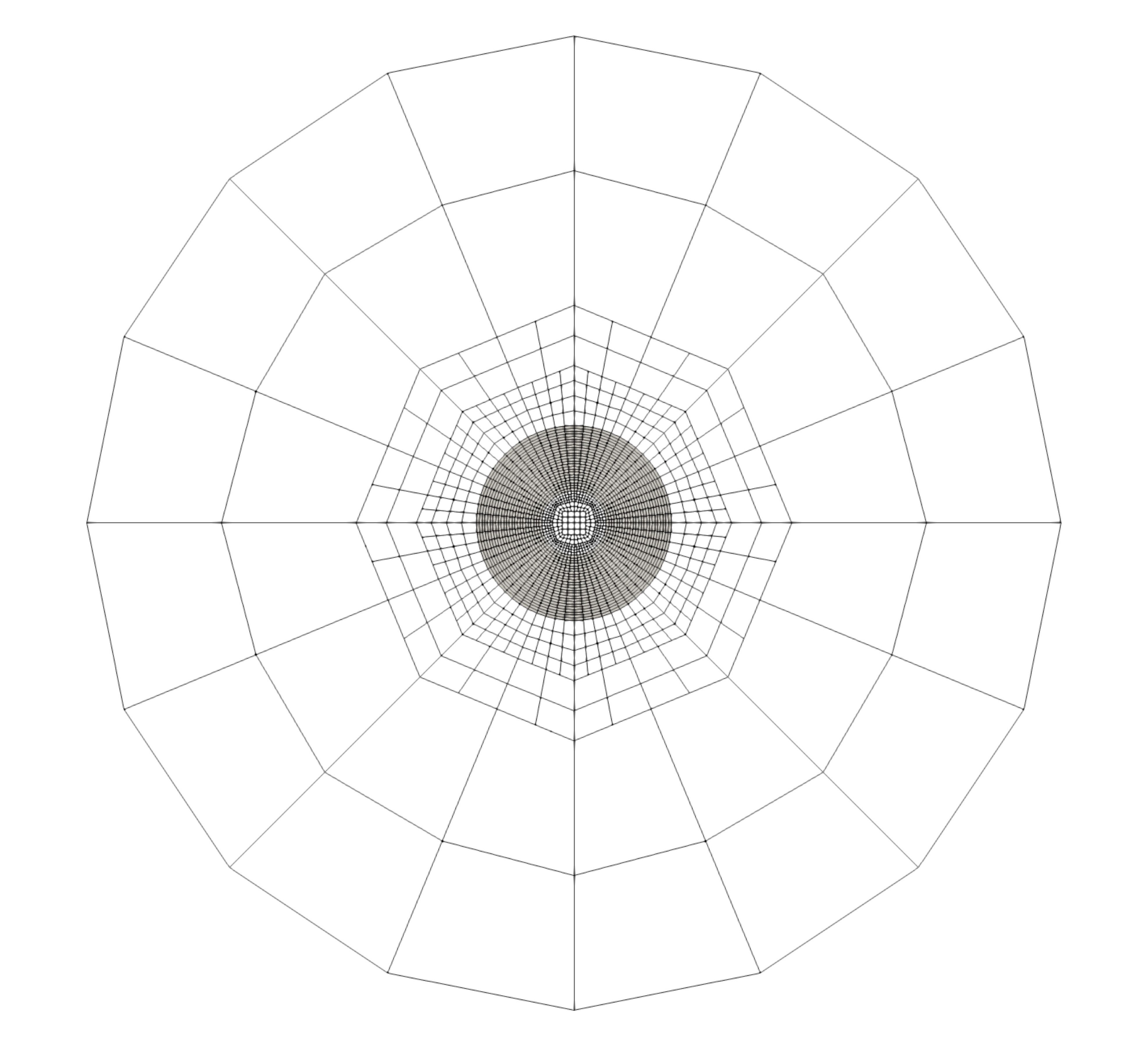} \\
				\includegraphics[scale=.09]{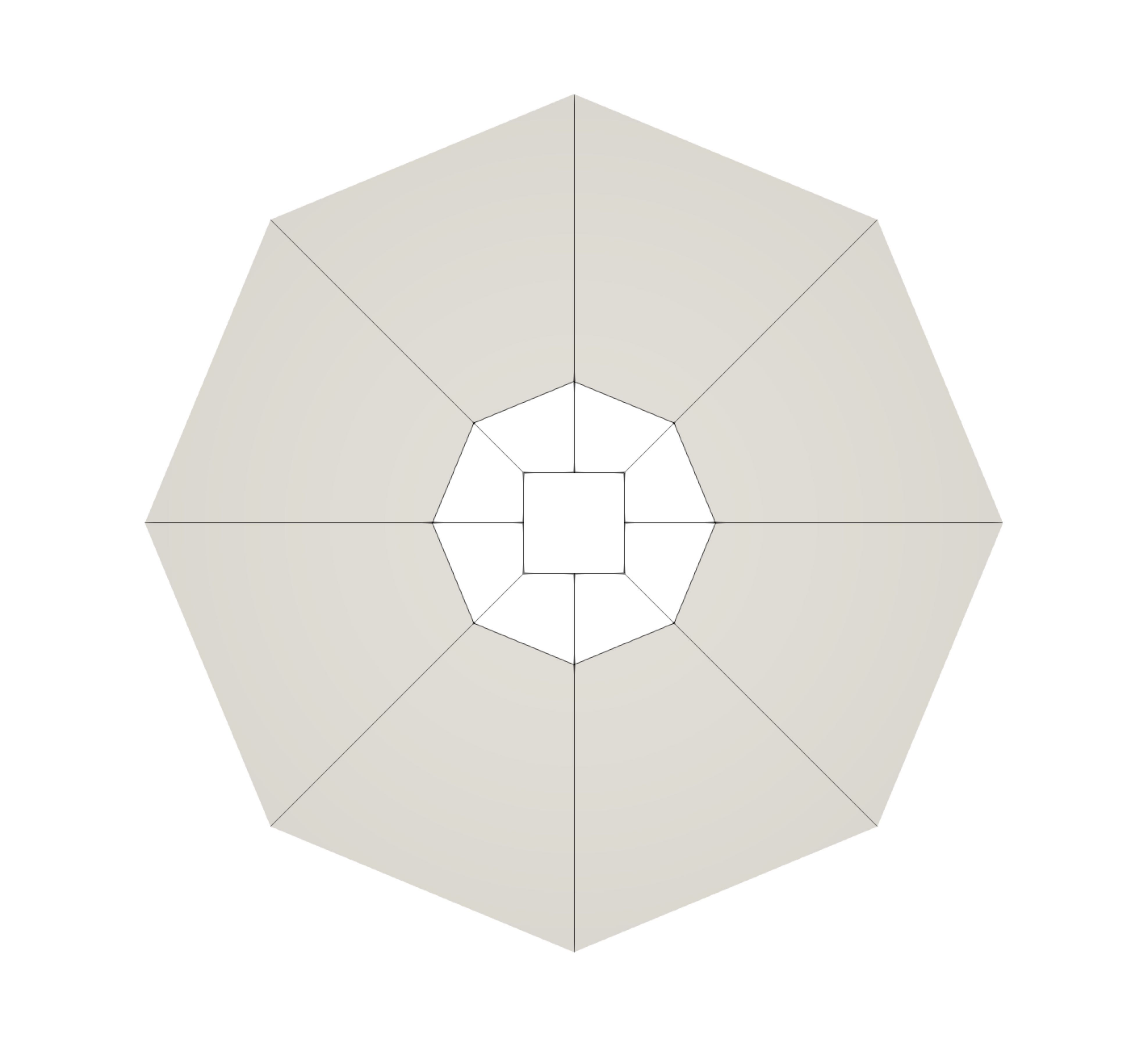} &
				\includegraphics[scale=.09]{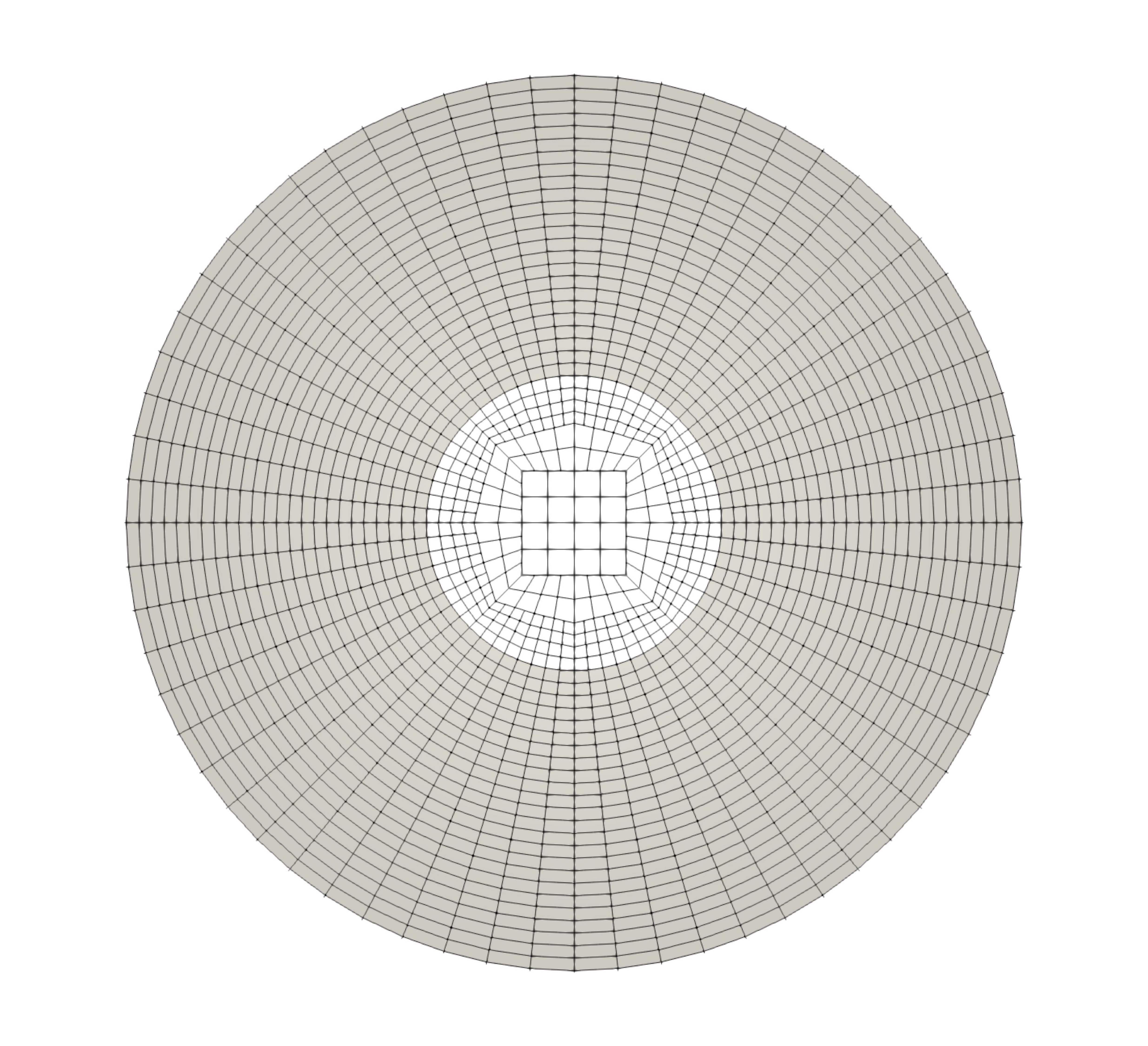} \\
				 (a) & (b)\\
			\end{tabular}
		\end{center}
		\caption{Polygonal approximation $\Omega_h^M(s)$ (first row) and zoom of the liquid region $\Omega_h$  (second row) with $M = 3$ and $s = 1$ associated to $\Omega$ with aspect ratio $\alpha = 0.33$. 
			The approximation of the liquid domain $\Omega_h$ is in gray.
			(a) Initial subdivision and (b) Three successive iterations of the refinement procedure.
		}
		\label{fig:mesh}
	\end{figure}

	 We end the discussion by describing the automatic generation of polygonal approximation $\Omega_h^M(s)$ of $\Omega^M(s)$ and its associated subdivision \modif{constrained to be an extension of a subdivisons of $\Omega$.}
	 We denote by $R^M(s)$ the radius of the truncated domain $\Omega^{M}(s)$.
	 We start with a coarse subdivision made of 25 quadrilaterals illustrated in \Cref{fig:mesh}(a).
	 To increase the accuracy while keeping the complexity under control, given $0<h\leq 1$, the refinement procedure consists of 
	\begin{itemize}
	\item refine the fluid domain $\Omega_h$ uniformly until all quadtilaterals have diameters smaller than $h$,  with possible additional refinements on $\Omega_h^M(s)$ to keep the number of hanging nodes to a maximum of one by edge.
	During this refinement process, newly created vertices at the boundary of $\Omega_h$ are placed on the boundary of $\Omega$.
	\item An exponential grading is performed outside $\Omega_h$,
	i.e. vertices on the azimuthal directions are placed at radii $r_j$ with
	$$
	r_j:=  e^{jh_0}/(1-\alpha), j = 1, 2, \cdots, \lceil M/h \rceil \quad \textrm{with} \quad h_0 = \ln(R^M(s)(1-\alpha))h/M;
	$$
	\end{itemize}
	We refer to  \Cref{fig:mesh}(b) for an illustration. 
	
	As in the case of infinite elecrodes configuration, the final resolution depends on the aspect ratio $\alpha$ to accommodate for the number of vortices during electroconvection.
	In the simulations presented in Section~\ref{s:simu}, we performed 5 uniform refinements for $\alpha \in \{0.1,0.2,0.33,0.452\}$, 6 for $\alpha\in \{0.56, 0.6446\}$ and $7$ for $\alpha=0.8$.
	
	\subsection{Approximation of the charge density~\eqref{eq:dl-EC-charge}}\label{subsec:charge}
	We now discuss the approximation of the surface charge density satisfying~\eqref{eq:dl-EC-charge}.
	Because $\Psi_K$ is harmonic in $\Omega$ in both electrodes configurations, 
	$$
	-\Delta_{\Omega} \Psi_\Omega 
	= -\Delta_{\Omega} \Psi_{0,\Omega}.
	$$
	Hence,  the surface charge density conservation relation reduces to
	$$
	\pt q + \bu \cdot\nabla_\Omega q - \Delta_\Omega \Psi_{0,\Omega} = 0 \qquad \textrm{in }\Omega \times \mathbb R_+.
	$$
	Notice that for a given electric potential $\Psi_{0,\Omega}$, the above partial differential equation for $q$ is a standard transport equation. 
	It is approximated with an explicit Runge-Kutta 2 method in time and standard finite elements in space. 
	
	 We recall $\mathcal T_h$ denotes a partition of a polygonal approximation $\Omega_h$ of the liquid domain $\Omega$ (see Section~\ref{subsec:pot}).
	 We introduce the finite element space for the surface charge density
	$$
	Q_h:= \{ v \in C^0(\overline{\Omega}_h) \ : \ v|_T \circ F_T \textrm{ is bilinear for all } T \in \mathcal T_h \textrm{ and }\int_{\Omega_h} v = 0 \}.
	$$
	Hence, given an approximation $\Phi_{0,h} \in V_h$ of $\Phi_{0,\Omega}$ determined as discussed in Section~\ref{subsec:pot}, the approximation $q_h(t) \in Q_h$ of $q(t)$ is defined as the solution to
	\begin{equation}\label{e:approx-charge-fem}
	\begin{split}
		&\int_{\Omega_h} \pt q_h(t) \varphi_h+ \int_{\Omega_h} \bu_h(t) \cdot\nabla_\Omega q_h(t)  \varphi_h\\
		& \qquad + \int_{\Omega_h} \nabla_\Omega\Psi_{0,h} \cdot  \nabla_\Omega \varphi_h
		- \int_{\partial \Omega_h} \nabla_\Omega \Psi_{0,h} \cdot \bnu_h \varphi_h
		= 0, \qquad \forall \varphi_h \in Q_h,
		\end{split}
\end{equation}
where $\bnu_h$ is the outside pointing normal to $\Omega_h$ (defined almost everywhere).
This relation is obtained upon multiplying the surface charge density equation and integrating by parts the electric potential term.
We suplement~\eqref{e:approx-charge-fem} with an approximated initial condition $\pi_h q(x,y,0)$, where $\pi_h$ denotes the $L^2(\Omega)$ projection onto $Q_h$.
	
	A Strongly Stability Preserving (SSP) two stage Runge-Kutta scheme proposed by~\cite{gottlieb2003SSP} is advocated for the time discretization.
	Set $q_h^0=\pi_h q(0)$ the $L^2$ projection of a given surface charge density onto $Q_h$ and $t_n:=n\tau$ for $n = 0,1, 2,\cdots$ for a time step parameter $\tau>0$. 
	Given approximation $q_h^n \in V_h$, $\Psi_{0,h}^n \in V_h$ (given by \eqref{e:final_Psi_infinite} or \eqref{e:final_Psi_slim} with $q$ replaced by $q_h^n$) and $\bu_h^n$ (see below) of the charge density $q(t_n)$, electric potential  $\Psi(t_n)$ and fluid velocity $\bu(t_n)$ respectively, 
	we approximate $q(t_{n+1})$ by $q_h^{n+1} \in Q_h$ as follows.
	In the first stage we seek $\mu_h^{(1)} \in Q_h$ satisfying
	\begin{equation}\label{e:approx-charge}
	\begin{split}
	\int_{\Omega_h} \mu_h^{(1)} \varphi_h
	=& \int_{\Omega_h} q_h^{n}  \varphi_h
	-\tau \int_{\Omega_h} \bu_h^{n}\cdot\nabla_\Omega q_h^n \varphi_h\\
	& - \tau \int_{\Omega_h} \nabla_\Omega \Psi_{0,h} \cdot \nabla_\Omega \varphi_h 
	+ \tau \int_{\partial \Omega_h} \nabla_\Omega \Psi_{0,h}\cdot \bnu_h \varphi_h,
	\end{split}
	\end{equation} 
	for all $\varphi_h\in Q_h$.
	In the second stage, we find $\mu_h^{(2)} \in Q_h$ solving~\eqref{e:approx-charge} but with $q_h^n$ replaced by $\mu_h^{(1)}$.		Then, we set
	\begin{equation}\label{e:q_n1}
	q_h^{n+1} := \frac{1}{2} (q_h^{n} +\mu_h^{(2)}).
	\end{equation}
	
	It is well documented that the finite element approximations of transport equation might be polluted by spurious oscillations.
	To circumvent this issue, we include the smoothness-based second-order maximum principle preserving viscosity method proposed in~\cite{Guermond2017invariant}. The numerical parameters used in the artificial viscosity are those recommended in (5.4) in~\cite{guermond2018sw}. 
	
	\subsection{Approximation of the Fluid Dynamic~\eqref{eq:dl-EC-momentum}} \label{subsec:ns}
	
	We discretize the fluid dynamic using backward differentiation scheme of order 2 (BDF-2) coupled with Taylor-Hood finite element approximations for the space discretization; see~\cite{Scott2007FEM}.
	Given a subdivision $\mathcal T_h$ of $\Omega_h$ constructed as in Section~\ref{subsec:pot},
	the finite element spaces for the velocity and pressure are defined by
	\begin{align*}
	\bW_h &:=  \{ \bv \in C^0(\overline{\Omega}_h)^2 \ : \  \bv|_T \circ F_T  \in \calQ_2^2,~ 
	\forall T\in \calT_h, ~ \bv|_{\partial \Omega_h} =0  \},\\
	X_h &:= \{ \theta\in C(\overline{\Omega}_h) 
	: \ \int_{\Omega_h} \theta = 0, \ \theta|_T\circ F_T \in \calQ_1, \forall T\in \calT_h \},
	\end{align*}
	where $\mathcal Q^i$, $i=1,2$, stands for the space of polynomial of (total) degree $i$.

	We start with $\bu_h^0 := \pi_h \bu(0)$, the $L^2$ projection of a given initial velocity $\bu(0)$ onto $\bW_h$.
	We assume that at time $t=t_n := n \tau$ (for the same time stepping parameter $\tau$ used for the surface charge density approximation), we have obtained $\Psi_h^{n} \in V_h$ given by \eqref{e:final_Psi_infinite} or \eqref{e:final_Psi_slim} with $q$ replaced by $q_h^{n}$ and $q_h^{n+1} \in Q_h$ given by \eqref{e:q_n1}.
	The  approximation  $(\bu_h^{n+1}, p_h^{n+1}) \in \bW_h \times X_h$ of $(\bu(t_{n+1}),p(t_{n+1}))$, $t_{n+1}:=(n+1)\tau$, is then defined as satisfying
	\begin{equation}
	\begin{split}
	  \int_{\Omega_h} \bu_h^{n+1} \cdot \bv_h 
	 &+ \frac{2}{3}\tau\calP \int_{\Omega_h} \nabla\bu_h^{n+1} \cdot \nabla\bv_h - \frac{2}{3}\tau \int_{\Omega_h} \DIV \bv_h  p_h^{n+1} \\
	 & = \int_{\Omega_h} \left(\frac{4}{3} \bu_h^{n} - \frac{1}{3} \bu_h^{n-1} \right) \cdot \bv_h
	- \frac{2}{3}\tau\int_{\Omega_h} \bv_h \cdot (\bu_h^n\cdot \nabla) \bu_h^{n} \\
	& \qquad \qquad \qquad -\frac{2}{3} \tau\calR\calP \int_{\Omega_h}  q_n^{n+1}\nabla_\Omega\Psi_h^{n} \cdot \bv_h,
	\end{split}
	\end{equation}
	and
	$$
	\int_{\Omega_h} \DIV \bu_h^{n+1} \theta_h = 0
	$$
	for all $(\bv_h, \theta_h) \in \bW_h\times X_h$. 		

	\section{Numerical Simulations}\label{s:simu}
	In this section we discuss the numerical results for two different settings with a particular emphasis on the effects of the Prandtl number
	$\mathcal P$, the Rayleigh number $\mathcal R$,  and the aspect ratio $\alpha$.

\begin{figure}[ht]
		\begin{center}
			\begin{tabular}{cc}
				\includegraphics[scale=.13]{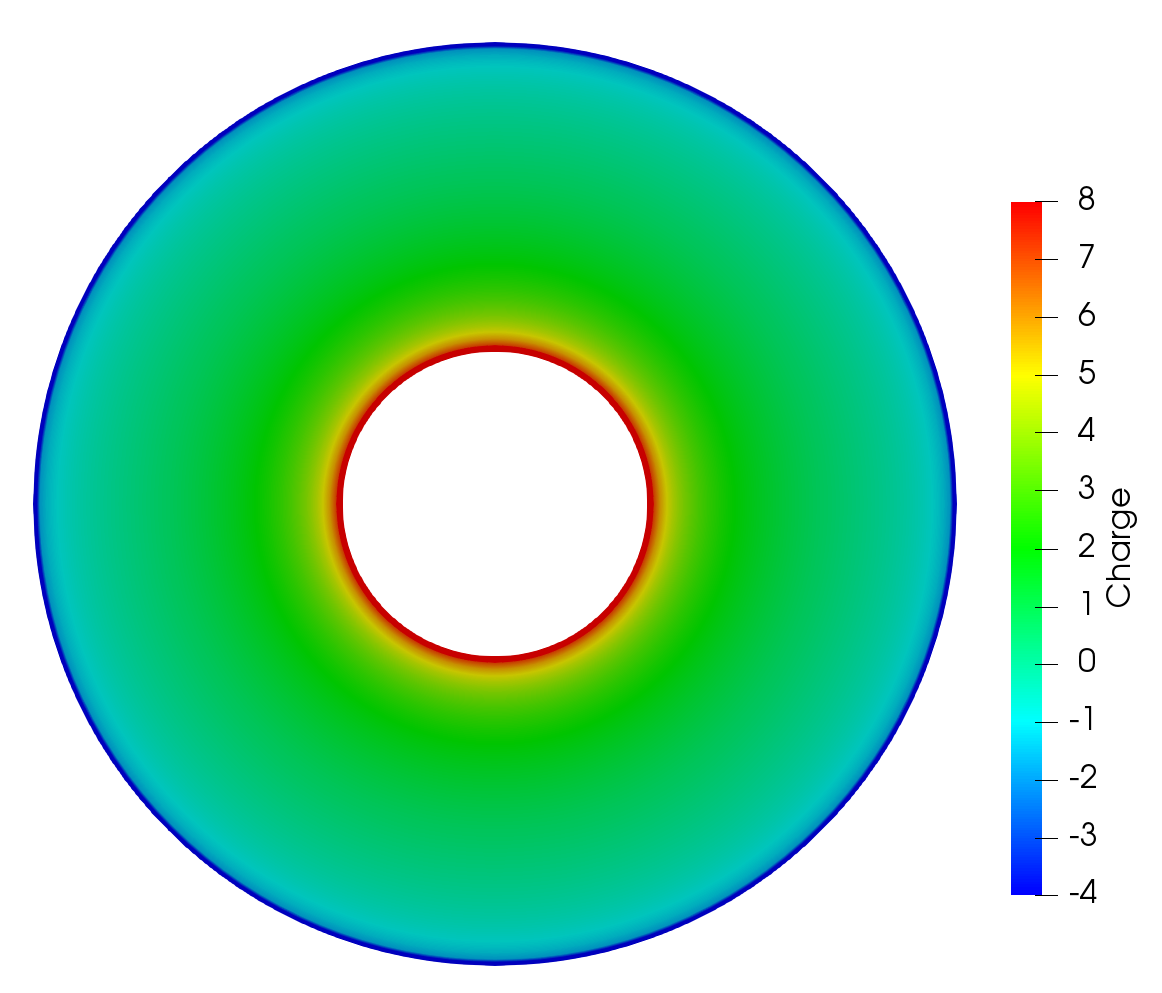}&
				\includegraphics[scale=.3]{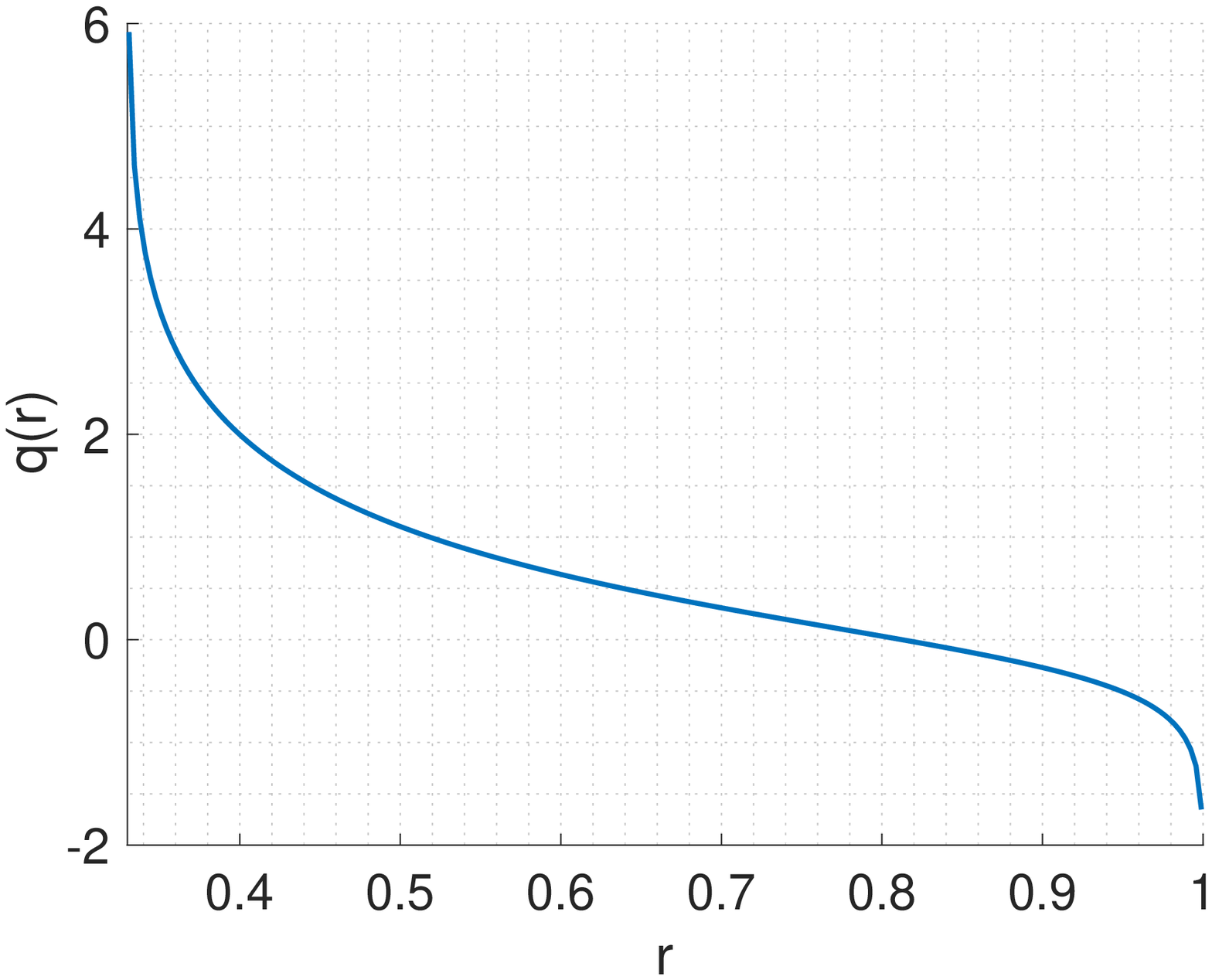} \\
				(a) & (b)
			\end{tabular}
		\end{center}
		\caption{The initial condition of charge density $q_0 = -g$ with $g$ defined by
			\eqref{d:g} for $\alpha = 0.33$. 
			(a) Initial surface charge density distribution;
			(b) Cross-section view of the initial surface charge density with $r=\sqrt{x^2+y^2}$.		
		}
		\label{fig:initial_q}
	\end{figure}

	The initial setting is common to all experiments. 
	The fluid is always starting at rest. \modif{For the surface charge density, we always start with an initial configuration in equilibrium with the electric potential associated with the slim electrode configuration. More precisely, the} initial charge density $q_0$ is chosen so that $\Psi = \Psi_K$ (i.e. $\Psi_{0,\Omega}=0$) in the decomposition~\eqref{e:decomp_slim}.
	This corresponds to setting $q_0 = -g$, where $g$ is given by ~\eqref{d:g}.
	We provide in Figure~\ref{fig:initial_q}  an illustration of the initial surface charge density for $\alpha = 0.33$.
	\modif{It is worth mentioning that we numerically observed the need of such compatibility between $q_0$ and $\Psi_0$ in the slim electrode setting but have not investigated this analytical question further.} 
	Notice that this configuration is unstable as the electric charges are aggregated near the outer boundary where the voltage is minimal.
	In the simulation below, we break the symmetry by adding a Gaussian white noise of magnitude no larger than $10^{-4}$.
	For comparison, we consider the same initial surface charge density in the infinite electrode case.
	
\begin{figure}[ht]
		\begin{center}
			\begin{tabular}{cc}
				\includegraphics[scale=.12]
				{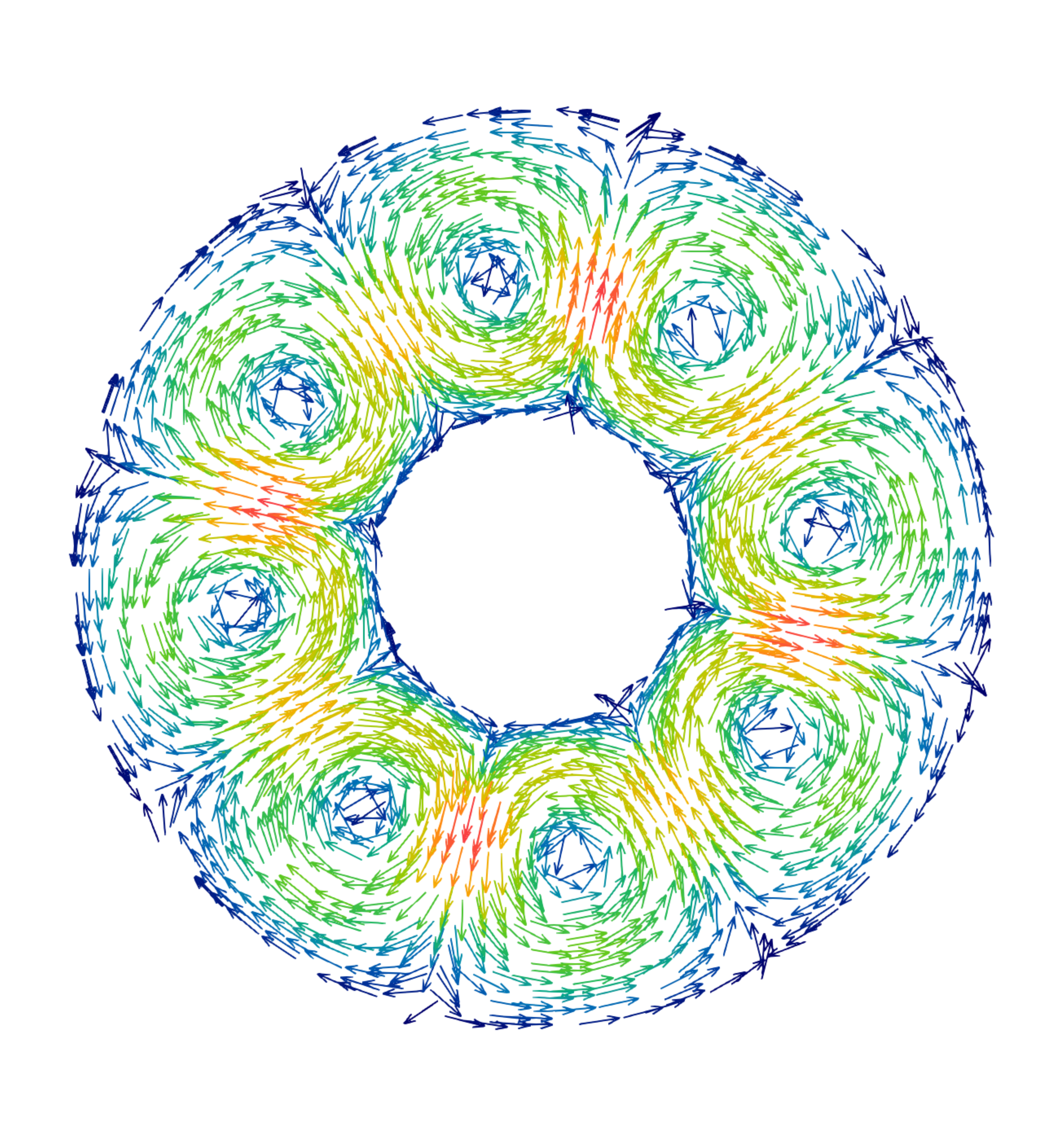}&
				\includegraphics[scale=.12]
				{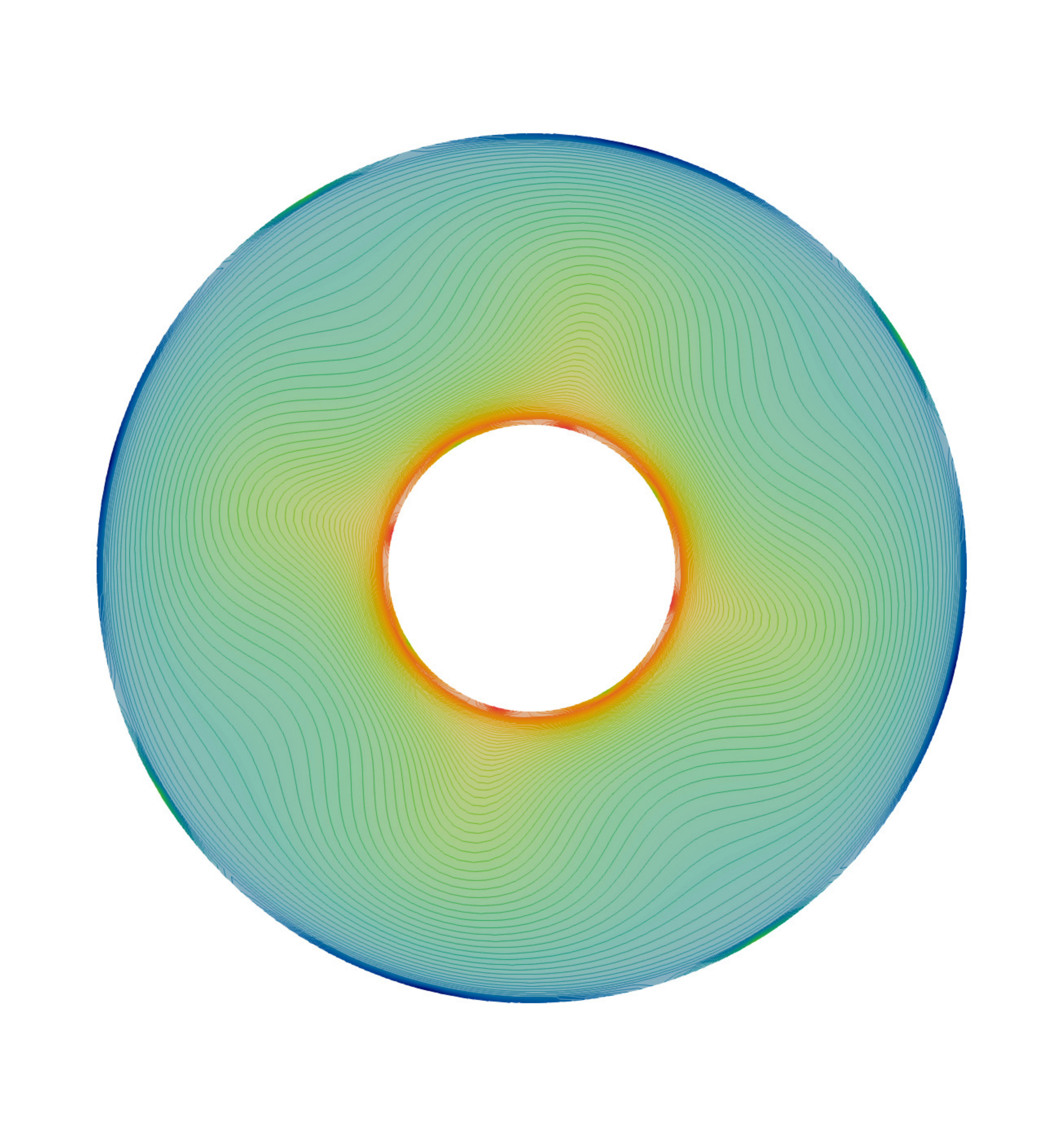} \\
				(a) & (b)
			\end{tabular}
		\end{center}
		\caption{Electroconvection for $\mathcal P = 10$, $\mathcal R = 100$, and $\alpha = 0.33$ at time $40$. 
			(a) Numerical approximation of the velocity field $\bu$;
			(b) Numerical approximation of the electric surface charge density distribution $q$.}
		\label{fig:alpha_0_33}
	\end{figure}

\begin{figure}[ht]
		\begin{center}
			\begin{tabular}{cc}
				\includegraphics[scale=.12]
				{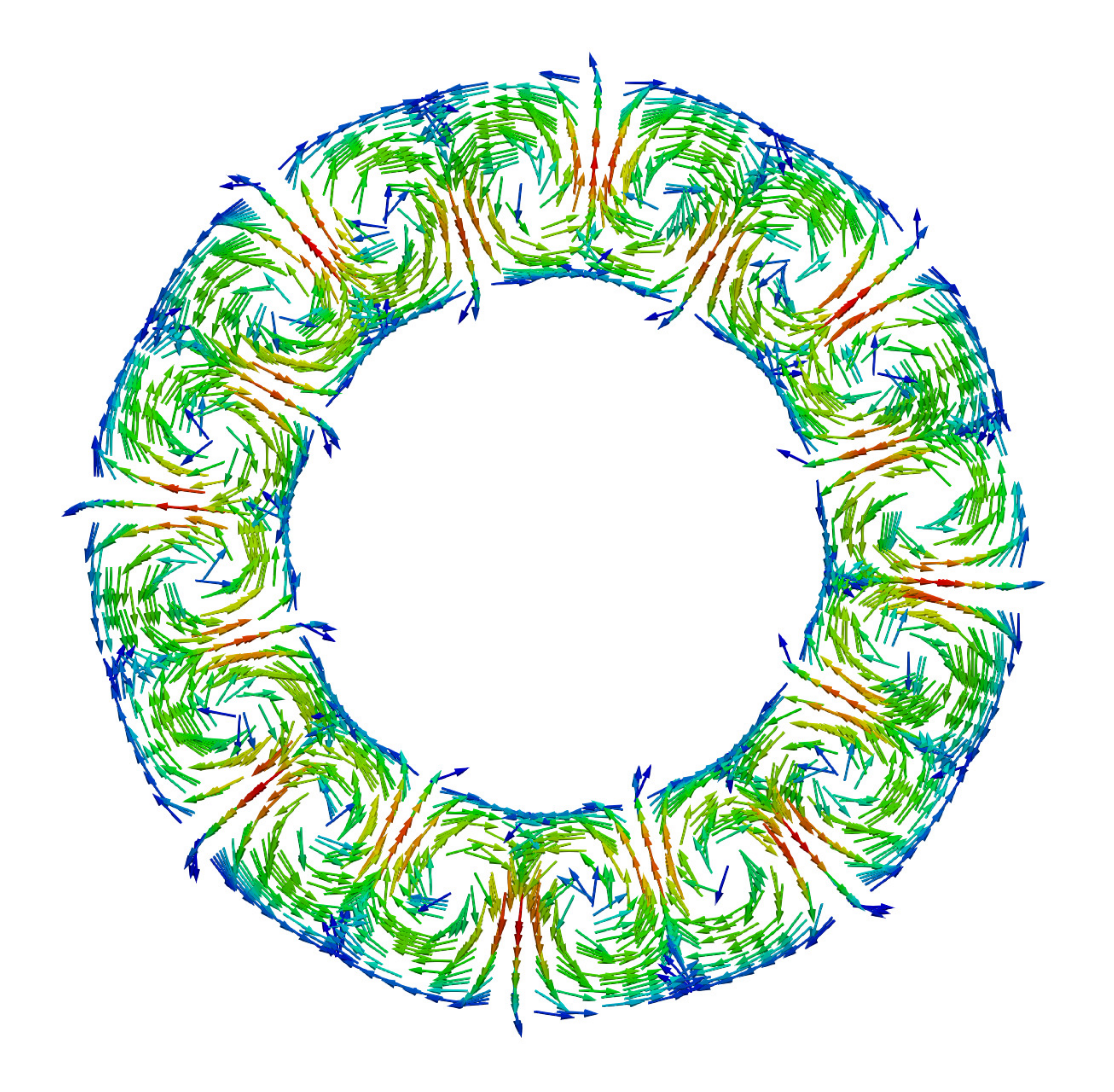}&
				\includegraphics[scale=.12]
				{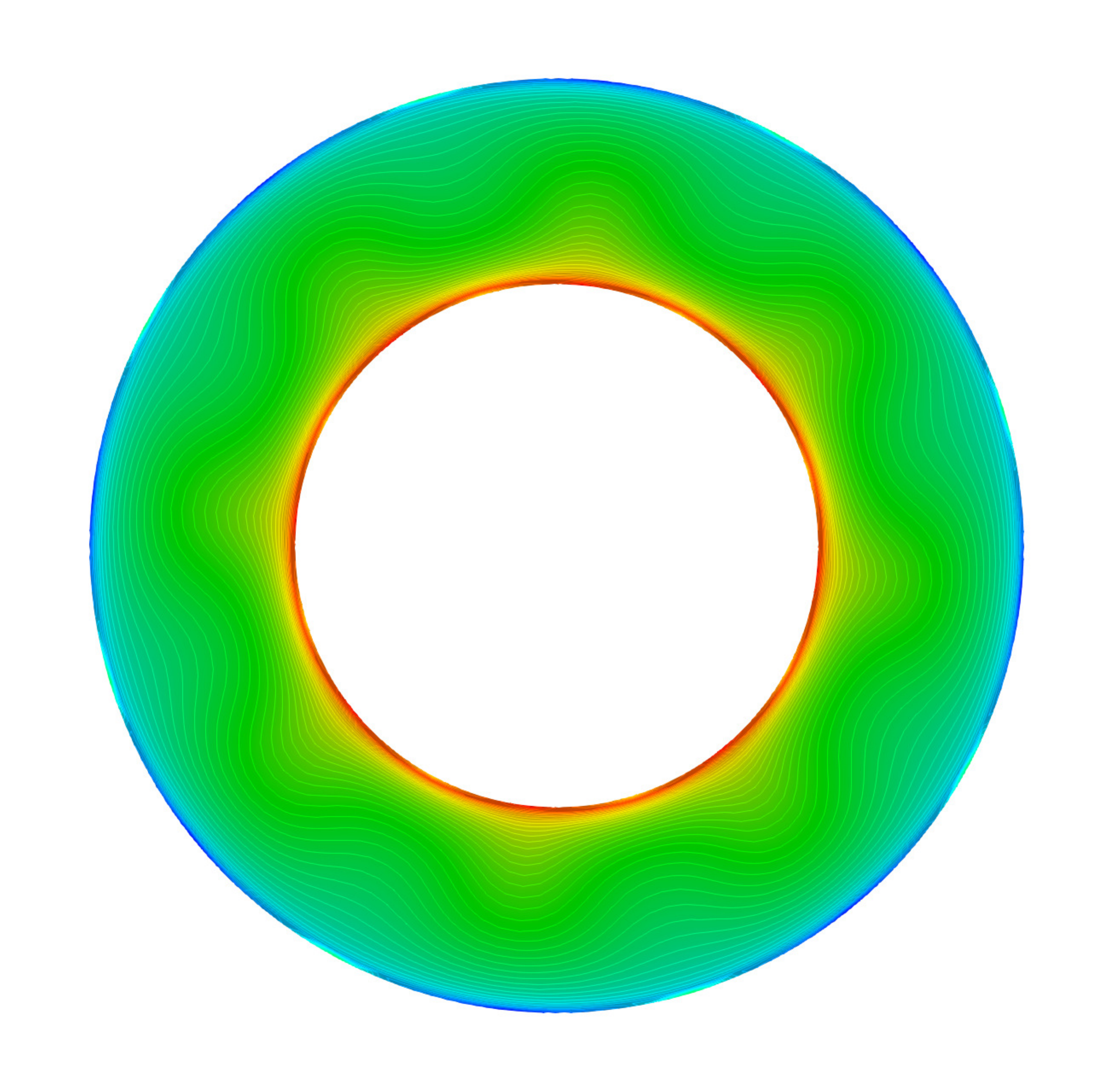} \\
				(a) & (b)
			\end{tabular}
		\end{center}
		\caption{Electroconvection for $\mathcal P = 10$, $\mathcal R = 100$, and $\alpha = 0.56$ at time $40$. 
			(a) Numerical approximation of the velocity field $\bu$;
			(b) Numerical approximation of the electric surface charge density distribution $q$.}
		\label{fig:alpha_0_56}
	\end{figure}

	In order to detect when and whether electroconvection occurs, we monitor two quantities: the kinetic energy
	\[
	E_k := \frac 1 2 \int_{\Omega} \abs{\bu}^2 \d\bx,
	\]
	and the circulation energy
	\[
	E_{curl} := \int_{\Omega} |\nabla\times \bu|^2\d\bx.
	\]
	Because of the added white noise to the initial surface charge density, the two quantities $E_k$ and $E_{curl}$ evolves initially.
	We declare the electroconvection phenomena to occur when both $E_k$ and $E_{curl}$ undergo a relative change greater than $0.1\%$ compared to their respective initial values before time $t=20$ (we ran our simulation further and did not observe any changes after that).  
	We define the critical Rayleigh number $\mathcal R_c$ as the threshold value for which fluids with lower Rayleigh numbers electroconvection do not occur (the Lorentz force  is not strong enough to overcome the electric and viscous dissipation).
	For Rayleigh numbers above $\calR_c$, the axis-symmetric distribution of charge density is broken and vortices appear. 
	The critical vortex pair number $\calN_c$ is the number of pairs of vortices during electroconvection.
	As an illustration, we provide in Figures~\ref{fig:alpha_0_33} and ~\ref{fig:alpha_0_56} the numerical approximation of sustained electroconvection with $\calN_c = 4$ and $\calN_c=8$ respectively.
	Similar structures but with different $\calN_c$ are observed for other different aspect ratios.
	
		We summarize the numerical parameters using in all the simulation below in Table~\ref{t:inf-parameter} for the infinite electrodes configuration and in Table~\ref{t:slim-parameter} for the slim electrodes configuration.
	\begin{table}[ht!]
				\begin{center}
		\begin{tabular}{|l|l|l|l|l|l|}
			\hline
			$\alpha$ & $\tau$ & DoFs for $\Psi_{0,\Omega}$ & DoFs for $q$ & DoFs for $\bu$ and $p$   \\ \hline
			0.1 & 0.001 & 4,224 & 4,224 & 33,280/4,224 \\ \hline
			0.2 & 0.001 & 5,280 & 5,280 & 41,600/5,280 \\ \hline
			0.33 & 0.001 & 7,392 & 7,392 & 58,240/7,392 \\ \hline
			0.452 & 0.001 & 9,504 & 9,504 & 74,880/9,504 \\ \hline
			0.56 & 0.001 & 12,672 & 12,672 & 99,840/12,672 \\ \hline
			0.6446 & 0.001 & 15,840 & 15,840 & 124,800/15,840 \\ \hline
			0.8 & 0.001 & 120,640 & 120,640 & 957,696/120,640 \\ \hline
		\end{tabular}
	\end{center}
	\caption{Parameter settings for infinite electrode simulations. Here ``DoFs'' stands for degrees of freedom.}
	\label{t:inf-parameter}
	\end{table}

	\begin{table}[ht!]
		\begin{tabular}{|l|l|l|l|l|l|}
			\hline
			$\alpha$ & $\tau$& $M$ & DoFs for $\Psi_{0,\Omega}$ & DoFs for $q$ & DoFs for $\bu$ and $p$   \\ \hline
			$\le 0.452$ & 0.001 & 3 & 9,009 & 6,272 & 49,664/6,272 \\ \hline
			0.56, 0.6446 & 0.001 & 3 & 24,960 & 24,960 & 198,144/24,960 \\ \hline
			0.8 & 0.001 & 3 & 99,072 & 99,072 & 789,504/99,072 \\ \hline
		\end{tabular}
	\caption{Parameter settings for slim electrode simulations. Here ``DoFs'' stands for degrees of freedom.}
	\label{t:slim-parameter}
	\end{table}

	\subsection{Comparison between the infinite and slim electrodes models}
\begin{figure}[h]
		\begin{center}
			\begin{tabular}{cc}
				\includegraphics[scale=.3] {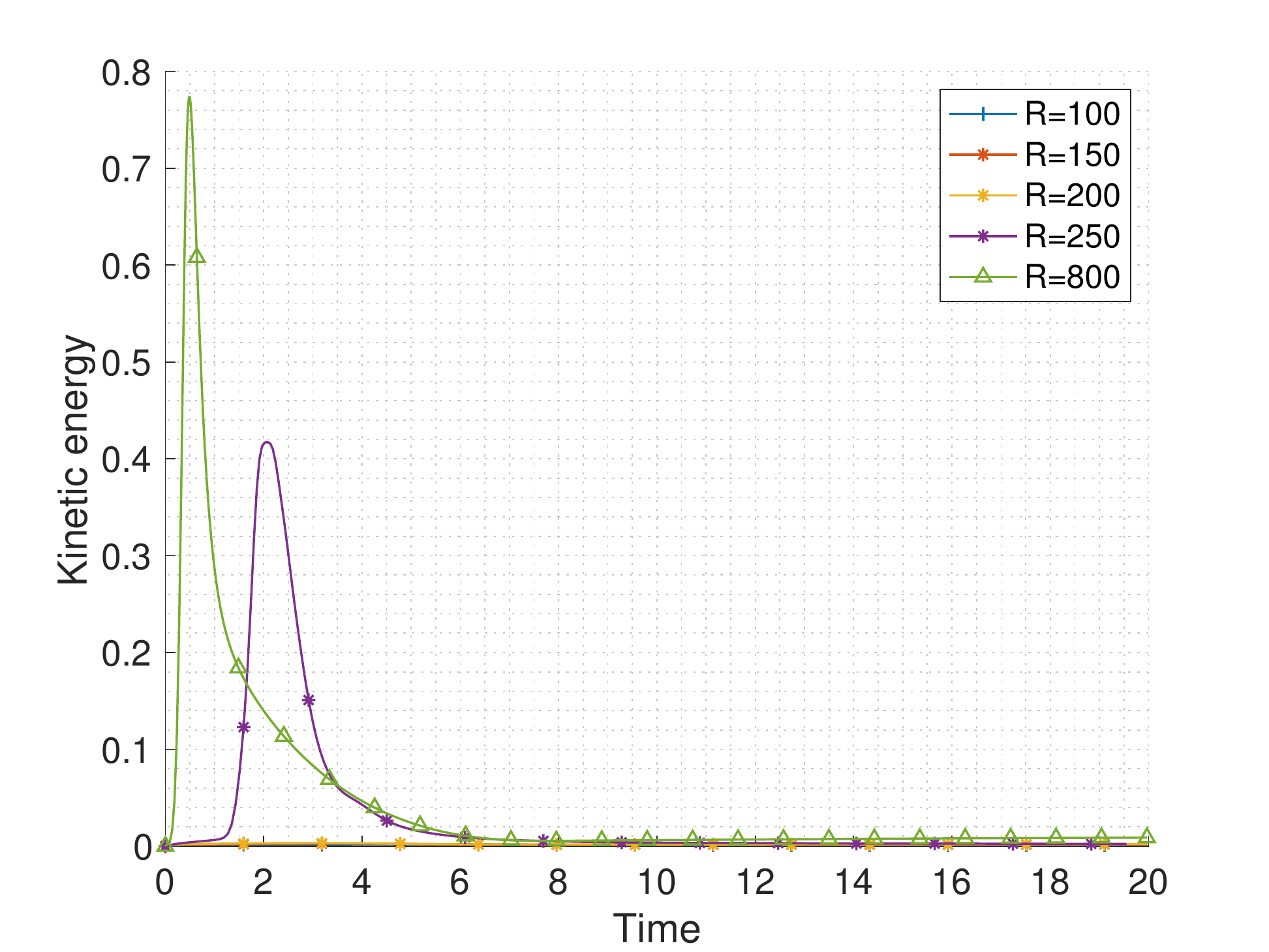}&
				\includegraphics[scale=.3] {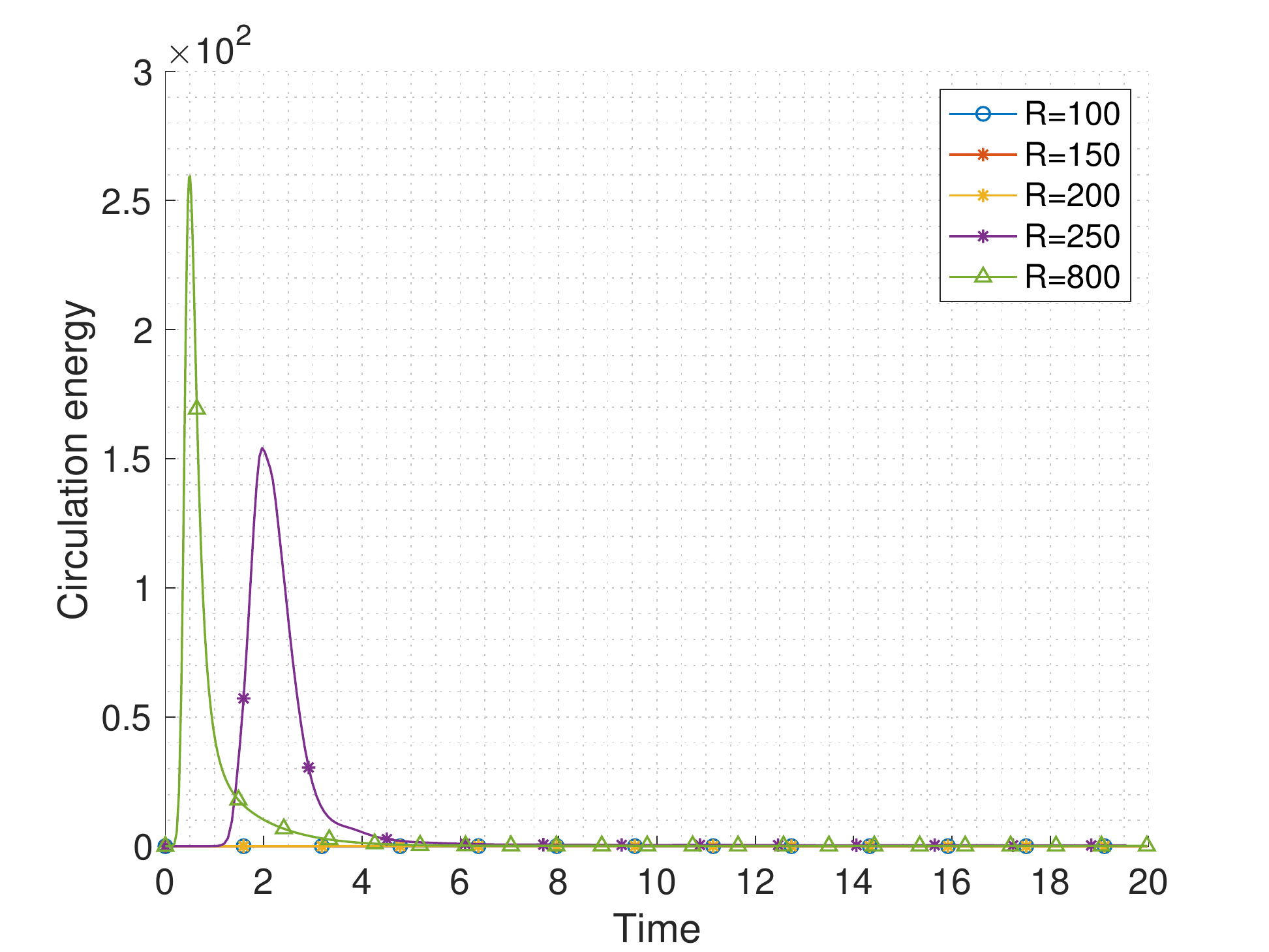} \\
				(a) & (b)
			\end{tabular}
		\end{center}
		\caption{Comparison of the Kinetic energy (a) and circulation enegry (b) for the infinite electrode model with $\mathcal P = 10$, $\alpha = 0.33$ and for several values of $\mathcal R$.}
		\label{pic:inf-R-033}
	\end{figure}

	For this comparison, we set the aspect ratio to $\alpha = 0.33$ and the Prandtl number to $\calP = 10$. 
	When the Rayleigh number is $\calR = 100$, electroconvection is observed in the slim electrodes configuration but not in the infinite electrodes configuration, which seems to requires more energy to trigger electroconvection. 
	In fact, electroconvection is observed in the infinite electrodes configuration for $\calR \ge 250$, see kinetic and circulation energies in~\Cref{pic:inf-R-033}.
	
\begin{figure}[ht]
		\begin{center}
			\begin{tabular}{cc}
				\includegraphics[scale=.3] {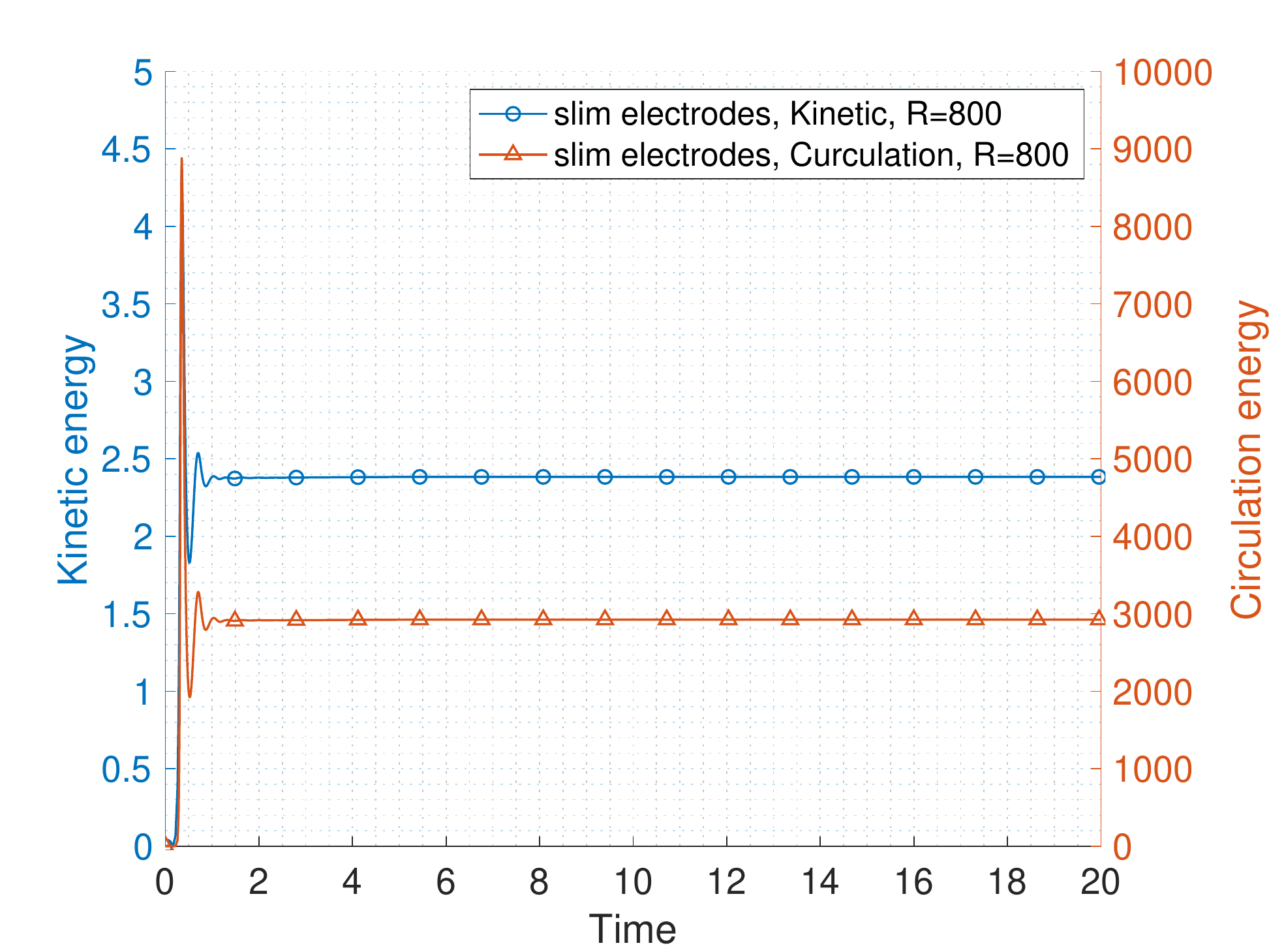}&
				\includegraphics[scale=.3] {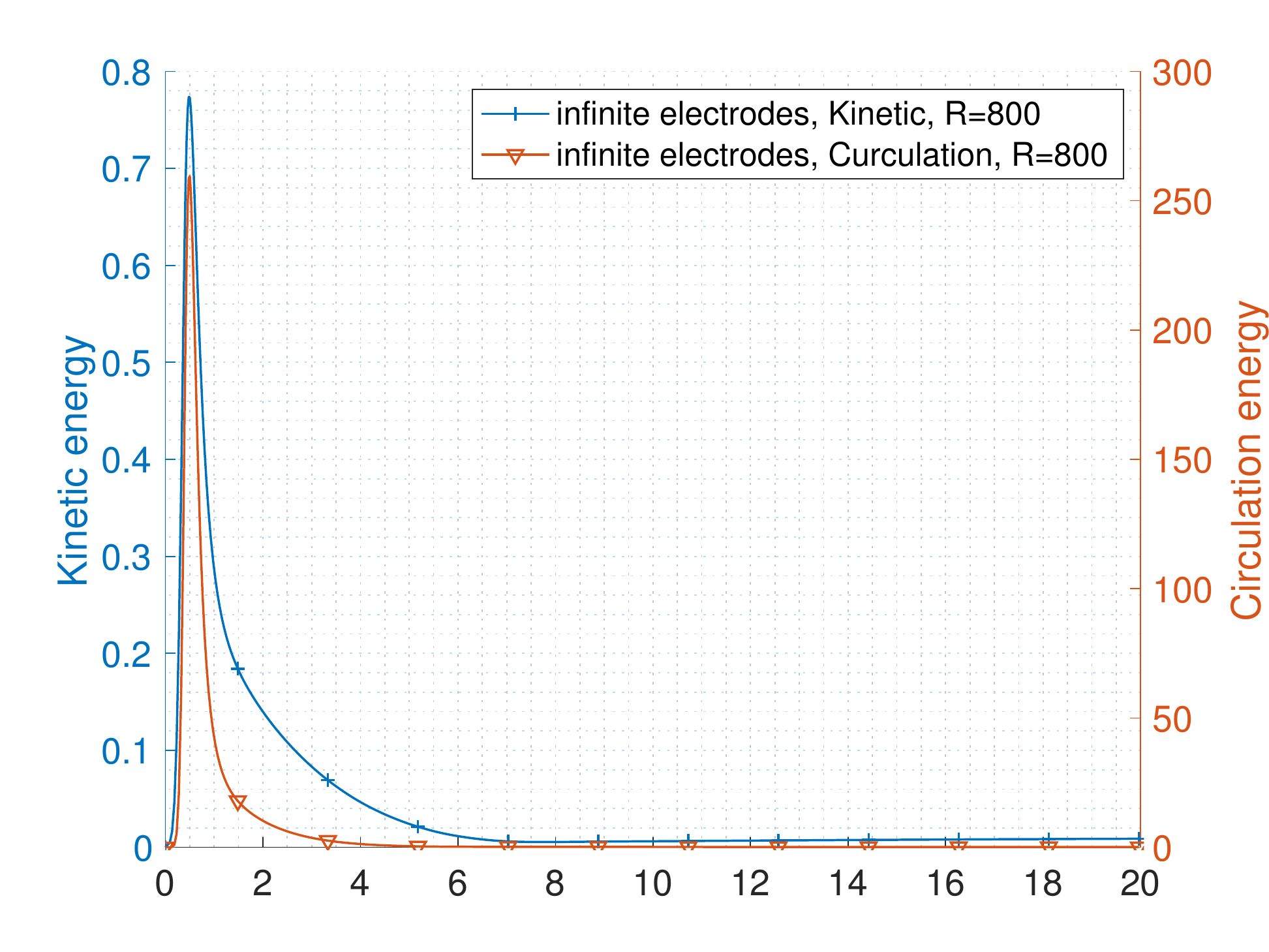} \\
				(a) & (b)
			\end{tabular}
		\end{center}
		\caption{Comparison of the evolution the energies $E_k$ and $E_{curl}$ versus time for the slim and infinite electrodes configurations with $\alpha = 0.33$, $\calR = 800$ and $\calP = 10$.
			The slim configuration energies are significantly larger and remain large during the entire evolution indicating a sustained electroconvection phenomena.}
		\label{fig:unbdd_bdd}
	\end{figure}

	Moreover, we find that even when electroconvection occurs in the infinite electrodes configuration, it cannot be sustained as in the slim electrode configuration. 
	To substantiate this fact, we set $\calP=10$ and $\calR=800$ and compare the energies in \Cref{fig:unbdd_bdd} for the two configurations.
	We observe that the energies in the slim case are not only significantly larger, they remain large as time evolves unlike in the infinite electrode configuration.
	
	The difference between the two models is striking. 
	In the presence of infinite electrodes model, electroconvection occurs as well but at much larger Rayleigh number $\calR$.
	Still, even when convective flows appears, they do not persist and quickly disappear, see Figure~\ref{fig:unbdd_bdd}.
	In contrast, the capability of maintaining a stable electroconvection phenomena indicates that the slim electrodes model is more adequate for electroconvection.
	Therefore, from now on we only consider the slim electrodes case.
	Worth mentioning, both models predict the same the number of pairs of vortices $\calN_c$, which seems indicating that $\calN_c$ depends mainly on the geometry as discussed later in Section~\ref{subsec:geometry}.	

	\subsection{Effect of the Prandtl number for the slim electrodes case}\label{s:effect_P}
\begin{figure}[h]
		\begin{center}
			\begin{tabular}{cc}
				\includegraphics[scale=.3] {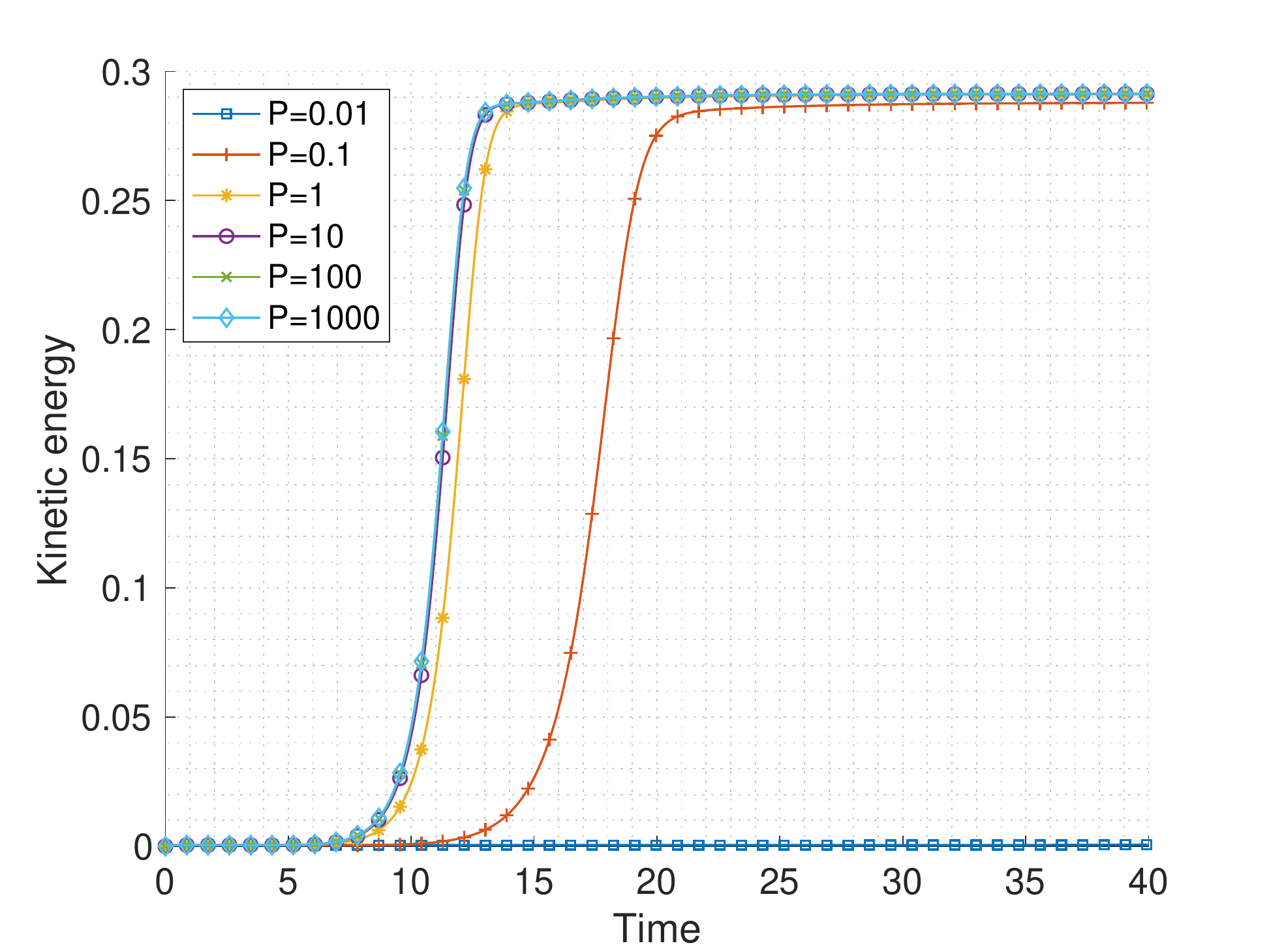}&
				\includegraphics[scale=.3] {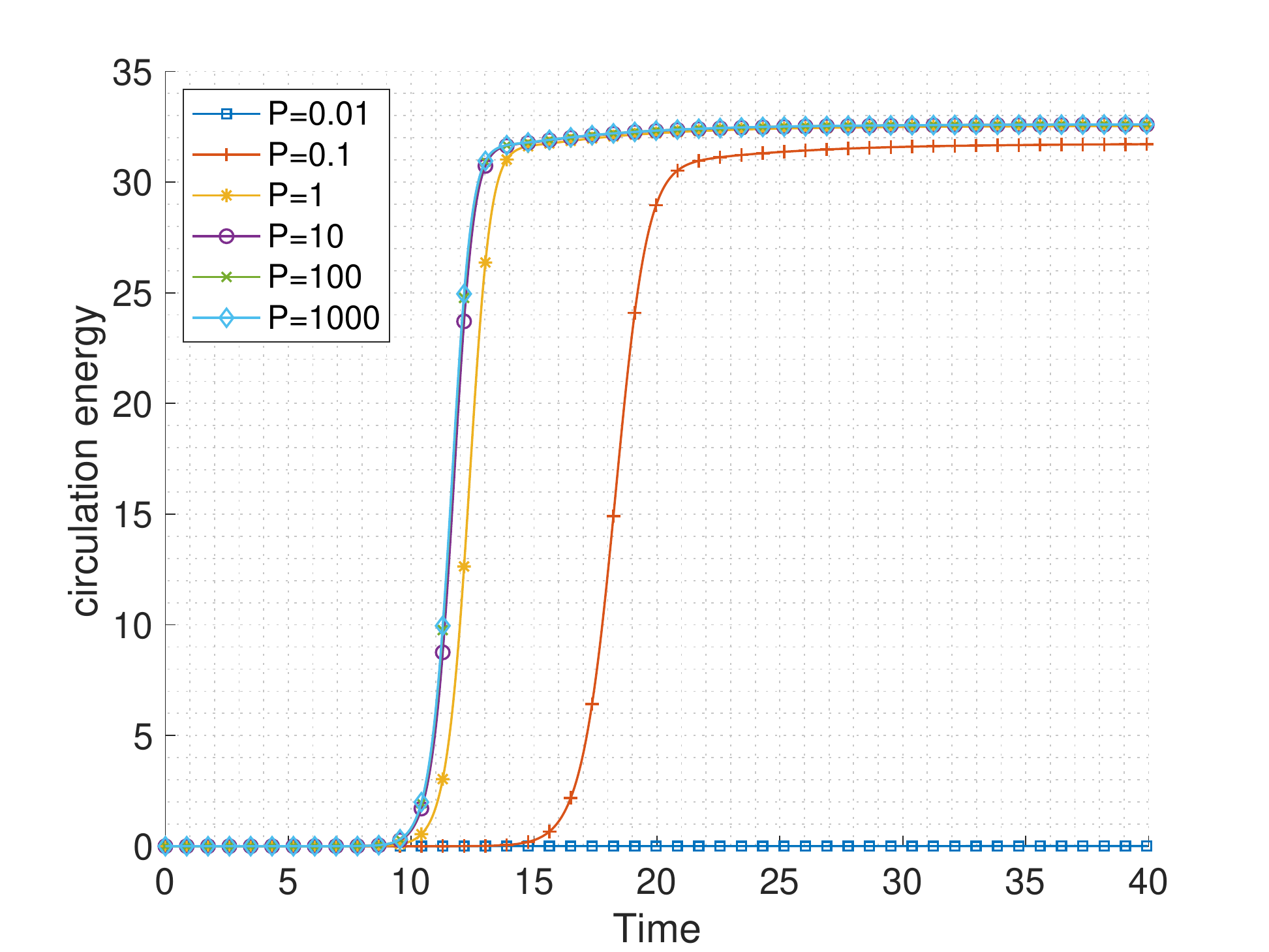} \\
				(a) & (b)
			\end{tabular}
		\end{center}
		\caption{Kinetic energy (a) and $E_{curl}$ (b) evolutions over time for different Prandtl number for $\alpha = 0.33$ and $\mathcal R = 100$.
			Increasing the Prandtl number decreases the time for the electroconvection to develop. 
			We observe that electroconvection cannot occur in this setting before $t=5$ irrespectively of the Prandtl number.
			Compare with~Figure~\ref{pics:R_033}.
		}
		\label{fig:P}
	\end{figure}

	The Prandtl number $\calP$ is the dimensionless ratio between the charge and viscous relaxation times.
	To understand its influence in the electroconvection phenomena, we fix $\alpha = 0.33$, $\mathcal R = 100$ and let $\mathcal P$ vary from \revtwo{$0.01$ to $1,000$.}
	In \Cref{fig:P} we report the kinetic energies and the circulation energies of the fluid.
	\revtwo{We observed that the electroconvection occurs for value of $\mathcal P \in[0.1,1000]$ which influences the activation time but not the long term behavior.}
	Worth noting, the energy plots in \Cref{fig:P} seem also to indicate that electroconvection cannot occur before $t=5$ irrespectively of the value of the Prandtl number.
	\revtwo{For smaller Prandtl number ($\calP = 0.01$ as shown in \Cref{fig:P}), electroconvection is not observed numerically. 
	The existence of a critical Prandtl number is in accordance with the nonlinear analysis provided in \cite{TDDW2007electroconvection}.
	For the rest of this section we will be focusing on the value $\calP = 10$ for the purpose of demonstrating convective flow.} 

	\subsection{Effect of the Geometry for the slim electrodes case}\label{subsec:geometry}
	The geometry of the annulus domain $\Omega$ is characterized by the aspect ratio $\alpha = R_i/R_o$; see~\eqref{e:ratio}.
	It turns out $\alpha$ affects the 	critical Rayleigh number $\mathcal R_c$ after which electroconvection occurs as well as the number of vortex pairs $\calN_c$.
	
\begin{figure}[ht]
		\begin{center}
			\begin{tabular}{cc}
				\includegraphics[scale=.3] {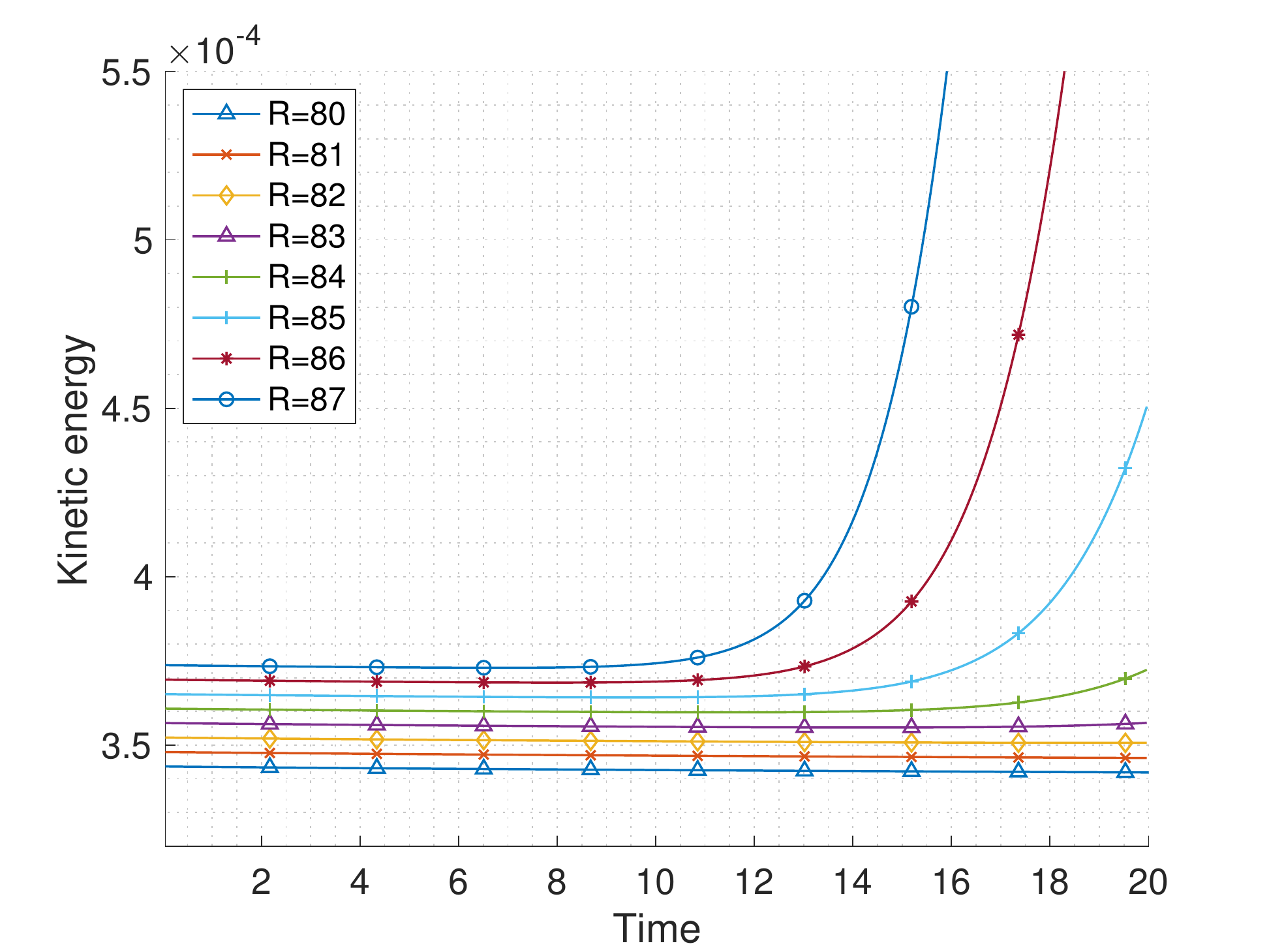}&
				\includegraphics[scale=.3] {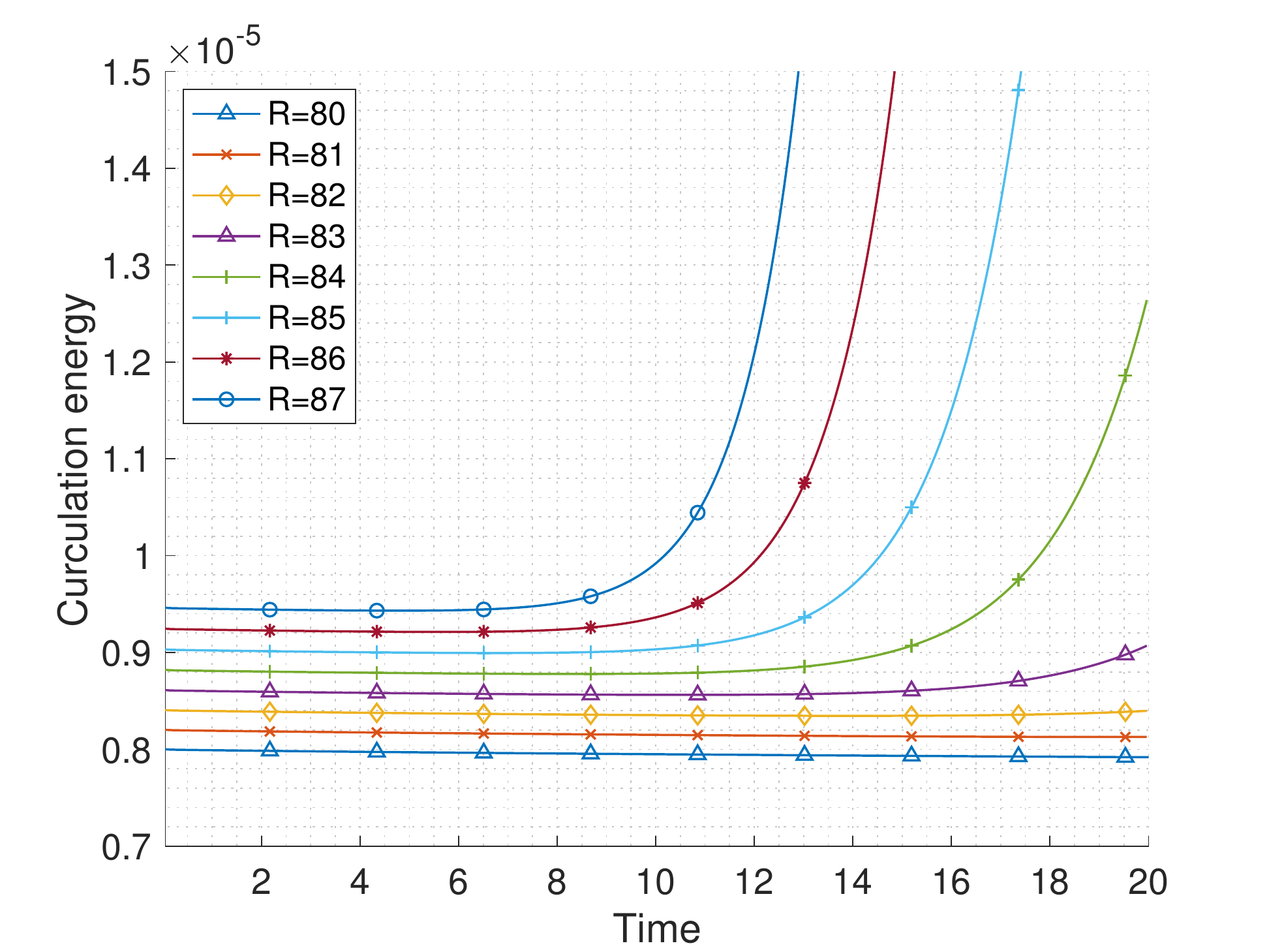} \\
				(a) & (b)
			\end{tabular}
		\end{center}
		\caption{Kinetic energy (a) and circulation energy (b) for slim electrodes with $\calP = 10$, $\alpha = 0.33$ and different values of $\calR $. }
		\label{pics:R}
	\end{figure}

	We fix $\calP = 10$ and determine for aspect ratio $\alpha \in [0.1,0.8]$ the corresponding critical Rayleigh number $\calR_c$.
	The latter is determined by running independent simulations starting from a relatively small $\calR$ at which electroconvection does not occur, and continuously increasing the value $\calR$ by one unit each time until electroconvection occurs.
	This allows us to determine the critical value $\calR_c$ up to $1$ unit.
	The system is considered steady when $E_k$ and $E_{curl}$ remains within $0.1\%$ relative difference throughout the simulation time.
	For example, \Cref{pics:R}(a) depicts the evolution of $E_k$ for 
	$\alpha = 0.33$ with $\mathcal R$ varying from $80$ to $87$. 
	The critical Rayleigh number satisfies $82 < \calR_c < 83$.
	
\begin{figure}[ht]
		\begin{center}
			\begin{tabular}{cc}
				\includegraphics[scale=.3] {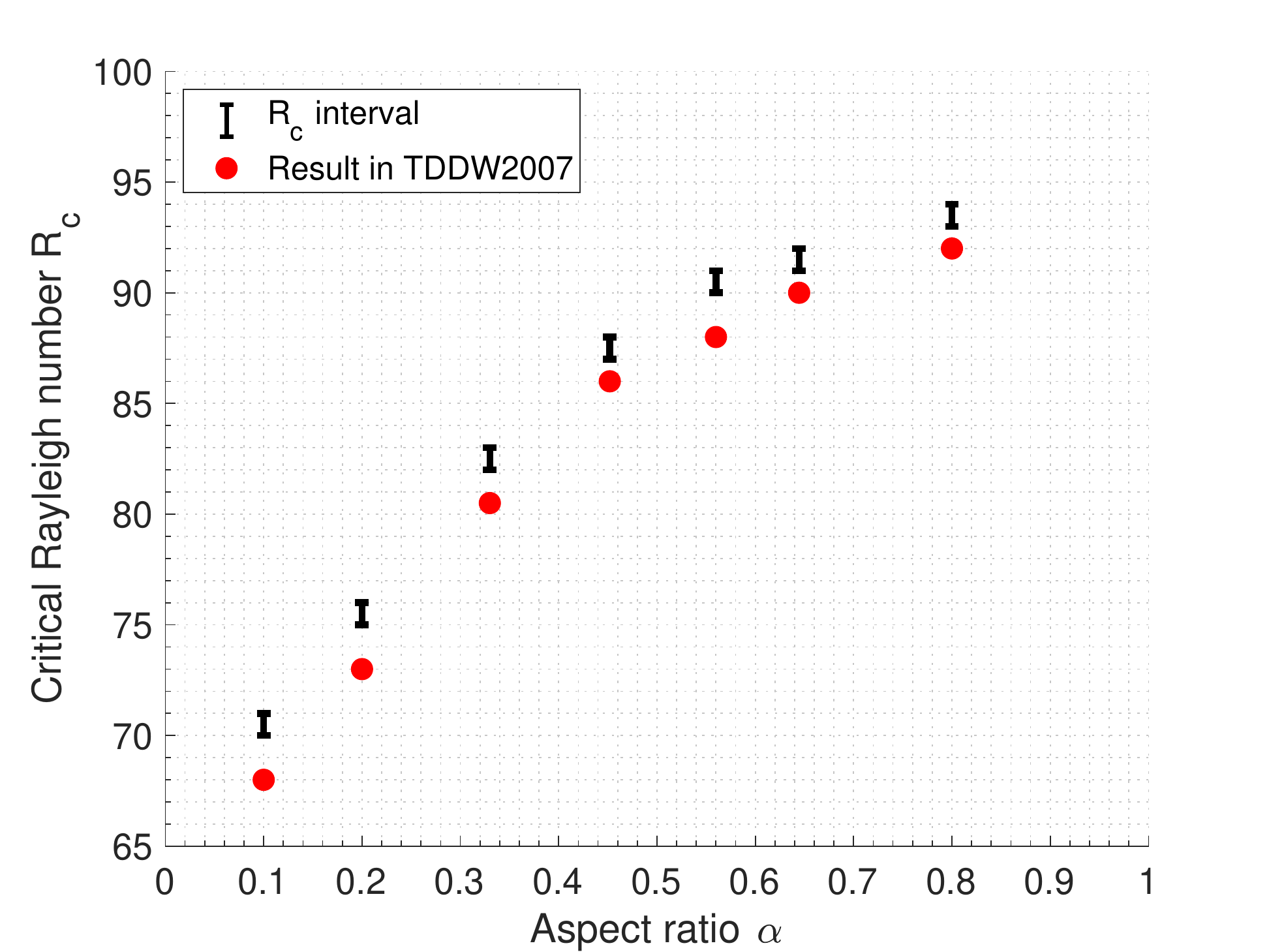}&
				\includegraphics[scale=.3] {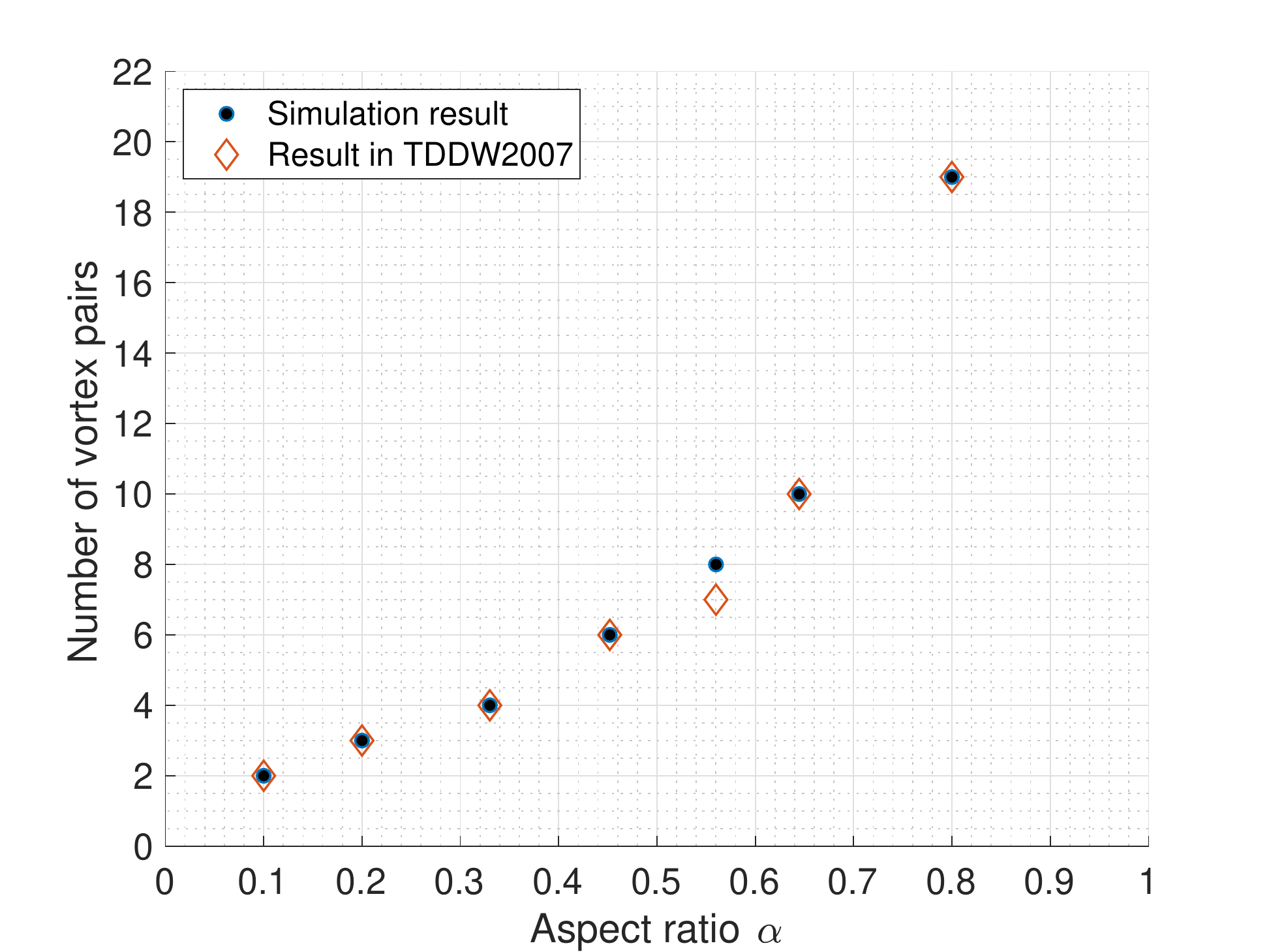} \\
				(a) & (b)
			\end{tabular}
		\end{center}
		\caption{Critical Rayleigh number (a) and number of vortex pairs (b) in the slim electrode case for $\mathcal P = 10$ and several aspect ratio $\alpha$. 
			Compare with the results provided in \cite{TDDW2007electroconvection}.
			The uncertainty intervals in part (a) are due to the increment used in the critical Rayleigh number exploration.
		}
		\label{pics:Rc}
	\end{figure}

	In \Cref{pics:Rc}(a) we report all the critical Rayleigh values for various $\alpha$ and compare them with  \cite{TDDW2007electroconvection}.
	\revtwo{Although not strictly matching the simulation results in Fig. 5 of \cite{TDDW2007electroconvection}, our results are in good agreement.
	In fact, in all cases, our numerical uncertainty intervals encompasses the theoretically predicted ranges reported in Fig. 5(a) of \cite{TDDW2007electroconvection}.  }
	
	The number of vortex pairs $\calN_c$ are strongly influenced by the geometry.
	To document this we fix $\mathcal P = 10$ and set, for each aspect ratio considered,  the Rayleigh number at the critical value $\mathcal R_c$.
	In \Cref{pics:Rc}(b) we compare the number of  vortex pairs $\calN_c$ obtained by our algorithm with results from
	\cite{TDDW2007electroconvection}.
	They match for all aspect ratio considered except for $\alpha = 0.56$ where they differ by one. 
	However, our predicted number of vortex pair in this case matches the experimental data Figure 5(b) of \cite{TDDW2007electroconvection}.

	\subsection{Effect of the Rayleigh number for the slim electrodes case} 
\begin{figure}[h]
		\begin{center}
			\begin{tabular}{cc}
				\includegraphics[scale=.3] {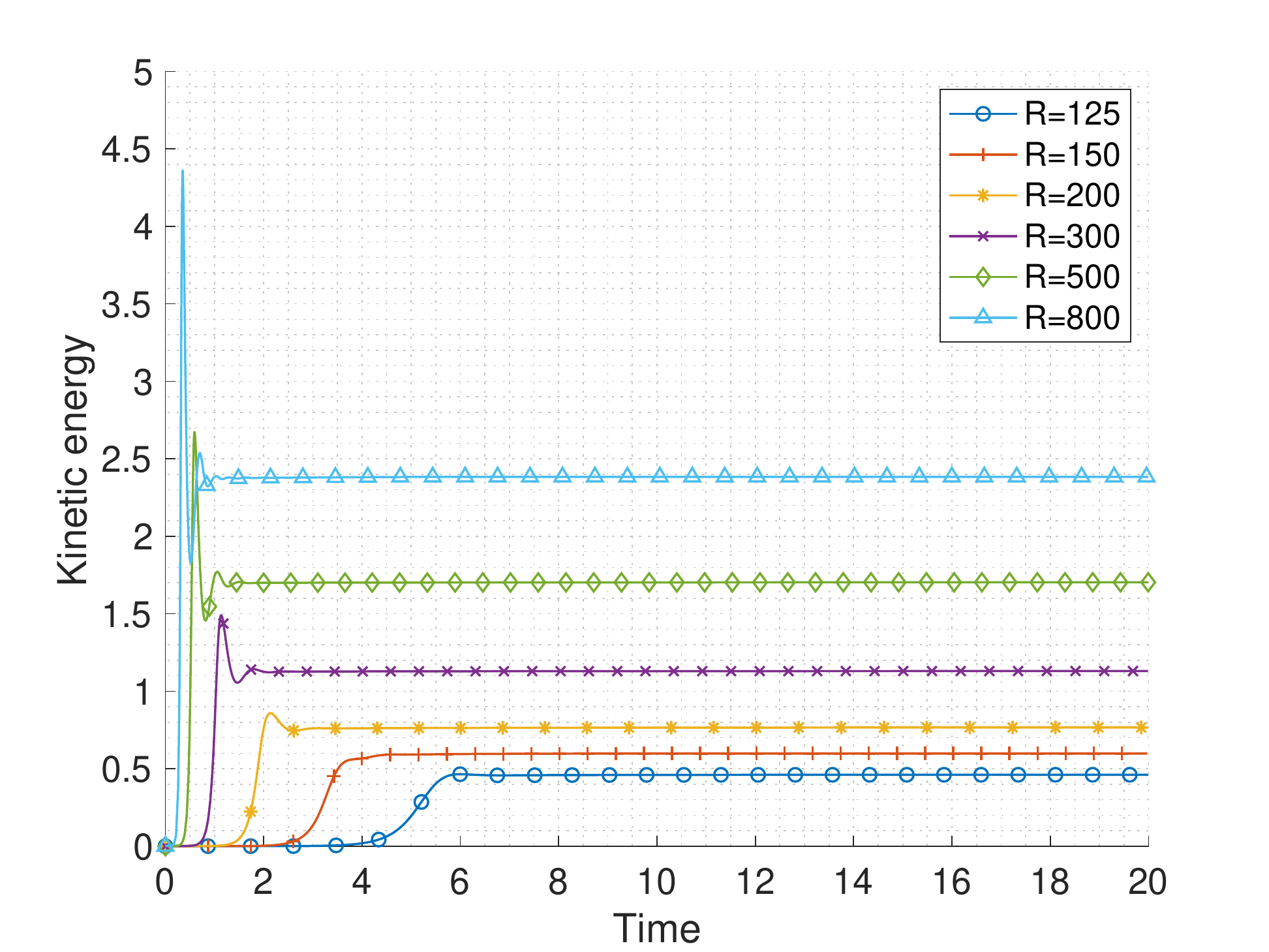}&
				\includegraphics[scale=.3] {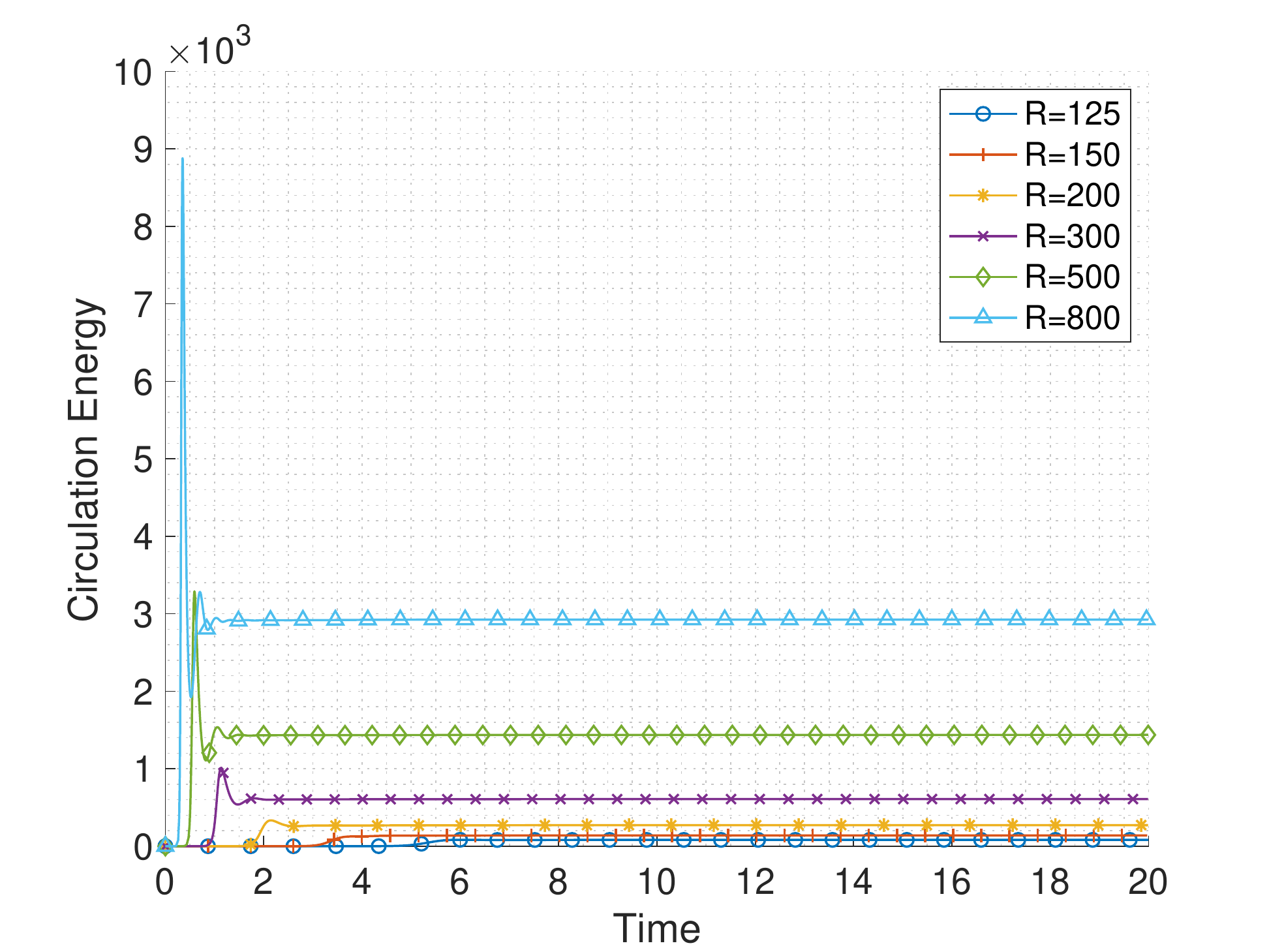} \\
				(a) & (b)
			\end{tabular}
		\end{center}
		\caption{Comparison of the Kinetic energy (a) and circulation energy (b) for the slim electrode model with $\mathcal P = 10$, $\alpha = 0.33$ and for several values of $\mathcal R$.
		}
		\label{pics:R_033}
	\end{figure}

	From the definition \eqref{PandR}, we realize that $\calR \propto V^2$.
	Increasing the Rayleigh number corresponds to a stronger electric field and thus a stronger Lorentz force. 
	In~\Cref{subsec:geometry},  we have already discussed the influence of the geometry on the critical Rayleigh value. 
	We now set $\alpha=0.33$ and $\calP=10$ and complete the investigation by increasing the value of $\mathcal R$ up to $10\mathcal R_c$. 
	The values of $E_k$ and $E_{curl}$ are reported in \Cref{pics:R_033}.
	We observe that larger Rayleigh numbers result in faster activation of the convection.
	This is similar as for the Prandtl number discussed in Section~\ref{s:effect_P} but in this case there does not seem to be a limiting time for before which electroconvection cannot occur. 
	Consequently fluids with larger Rayleigh numbers  develop a stable electroconvection at earlier times.
	The strength of convection is also stronger for larger Rayleigh number.
	
	\section{Conclusions} \label{s:conclusion}

In this paper, we have derived a mathematical model for electrically driven convection in an annular two dimensional fluid.
The two different electrodes configurations are considered: infinite and slim. 
Depending on the electrodes configuration, nonlocal representations of the electric potential are derived on the liquid domain.
This together with the surface density charge conservation relation and the Navier-Stokes system for the fluid dynamics yield a system of partial differential equations defined only in the two dimensional and bounded liquid region.

The proposed numerical methods take advantage of this dimension reduction and only require  discretization tools readily available in most finite element codes. 
Our numerical simulations reveal that the slim electrodes configuration is more favorable for electroconvection: it requires less energy and is able to sustain the effect. 
It is therefore the configuration chosen to provide a numerical study of the three nondimensional parameters describing the electroconvection system.

\revtwo{We observe the Prandtl number must be above a critical value for electroconvection occurs.
However, values above such critical number only affects the electroconvection activation time but not the long term behavior.}
\modif{For several geometries, we also exhibit  critical Rayleigh numbers above which electroconvection occurs.}

\bibliographystyle{abbrvnat} 
\bibliography{ref}
	
\end{document}